\documentclass[11pt]{article}
\usepackage{amsmath,amssymb,amsfonts}
\usepackage[latin1]{inputenc}
\usepackage{color}
\numberwithin{equation}{section}\newtheorem{thm}{Theorem}[section]

\newtheorem{alem}[thm]{Lemma}
\newtheorem{aprop}[thm]{Proposition}
\newtheorem{acor}[thm]{Corollary}
\newtheorem{arem}[thm]{Remark}

\newenvironment{dem}[1][]%
   {\ \\ {\bf Proof#1. }}%
   {\hfill\mbox{\rule{2 true mm}{3 true mm}}}
\newenvironment{adem}[1][]%
   {\ \\ {\bf Proof #1. }}%
   {\hfill\mbox{\rule{2 true mm}{3 true mm}}}
   {\ \\ {\bf Example #1. }}%
   {\hfill\mbox{\rule{2 true mm}{3 true mm}}}
\newcommand{\R}{\mathbb{R}}

\newcommand{\E}{{\mathbb E}}

\newcommand{\N}{{\mathbb N}}

\newcommand{\mycomments}[1]{  }

\newcommand{\be}{ \begin{eqnarray*}  }
\newcommand{\ee}{ \end{eqnarray*}  }

\DeclareMathOperator{\Tr}{Tr}

 \title{Optimal transport bounds between the time-marginals of a multidimensional diffusion and its Euler scheme}
\author{A. Alfonsi, B. Jourdain\thanks{Universit\'e Paris-Est, CERMICS, Projet MathRisk
    ENPC-INRIA-UMLV, 6 et 8 avenue Blaise Pascal, 77455 Marne La Vall\'ee, Cedex
    2, France, e-mails : alfonsi@cermics.enpc.fr, jourdain@cermics.enpc.fr. This research benefited
    from the support of the ``Chaire Risques Financiers'', Fondation du
    Risque, the French National Research Agency (ANR) under the program
 ANR-Stab and the Labex Bézout.
}~~and A. Kohatsu-Higa \thanks{Ritsumeikan University and Japan Science and Technology Agency, Department of Mathematical Sciences, 1-1-1
Nojihigashi, Kusatsu, Shiga, 525-8577, Japan. e-mail: arturokohatsu@gmail.com. This research was supported by grants of the Japanese goverment.}}

\begin{document}
\maketitle\begin{abstract}
In this paper, we prove that the time supremum of the Wasserstein distance between the time-marginals of a uniformly elliptic multidimensional diffusion with coefficients bounded together with their derivatives up to the order $2$ in the spatial variables and H\"older continuous with exponent $\gamma$ with respect to the time variable and its Euler scheme with $N$ uniform time-steps is smaller than $C \left(1+\mathbf{1}_{\gamma=1} \sqrt{\ln(N)}\right)N^{-\gamma}$. To do so, we use the theory of optimal transport. More precisely, we investigate how to apply the theory by Ambrosio {\it et al.} \cite{ags} to compute the time derivative of the Wasserstein distance between the time-marginals. We deduce a stability inequality for the Wasserstein distance which finally leads to the desired estimation.
\end{abstract}
\section{Introduction}
Consider the $\R^d$-valued Stochastic Differential Equation (SDE) :
\begin{equation}
X_t=x_0+ \int_0^t b(s,X_s)ds + \int_0^t \sigma(s,X_s)dW_s,\;t\leq T\label{sde}
\end{equation}
with $T>0$ a finite time-horizon, $(W_t)_{t\in[0,T]}$ a $d$-dimensional standard Brownian motion, $b:[0,T]\times\R^d\to \R^d$ and $\sigma : [0,T]\times\R^d\to {\cal M}_d(\R)$ where ${\cal M}_d(\R)$  denotes the set of real $d\times d$-matrices. In what follows, $\sigma$ and $b$ will be assumed to be Lispchitz continuous in the spatial variable uniformly for $t\in[0,T]$ and such that $\sup_{t\in[0,T]}(|\sigma(t,0)|+|b(t,0)|)<+\infty$ so that trajectorial existence and uniqueness hold for this SDE. 

We now introduce the Euler scheme. To do so, we consider for $N\in \N^*$ the regular time grid $t_i=\frac{iT}{N}$. We define the continuous time Euler scheme by
the following induction for $i\in\{0,\hdots,N-1\}$ :
\begin{equation}\label{schemeul}
\bar{X}_{t}=\bar{X}_{t_i}+b(t_i,\bar{X}_{t_i})(t-t_i)+\sigma(t_i,\bar{X}_{t_i})(W_t-W_{t_i}), \ t\in[t_i,t_{i+1}],
\end{equation}
with $\bar{X}_{t_0}=x_0$.
By setting $\tau_t=\lfloor \frac{Nt}{T}\rfloor\frac{T}{N}$, we can also write the Euler scheme as an It\^o process :
\begin{equation}\label{eul_conti} \bar{X}_{t}=x_0+ \int_0^t b(\tau_s,\bar{X}_{\tau_s})ds + \int_0^t \sigma(\tau_s,\bar{X}_{\tau_s})dW_s,\;t\leq T. 
\end{equation}

The goal of this paper is to study the Wasserstein distance between the laws ${\cal L}(X_t)$ and ${\cal     L}(\bar{X}_t)$ of $X_t$ and $\bar{X}_t$. We first recall the definition of the Wasserstein distance. Let $\mu$ and $\nu$ denote two probability measures on~$\R^d$ and $\rho\ge 1$. The $\rho$-Wasserstein distance between $\mu$ and $\nu$ is defined by
\begin{equation}
   {\cal W}_{\rho}(\mu,\nu)=\left(\inf_{\pi\in\Pi(\mu,\nu)}\int_{\R^d\times\R^d}|x-y|^\rho\pi(dx,dy)\right)^{1/\rho},\label{defwas}
\end{equation}
where $\Pi(\mu,\nu)$ is the set of probability measures on $\R^d\times\R^d$ with respective marginals $\mu$ and $\nu$. In this paper, we will work with the Euclidean norm on~$\R^d$, i.e. $|x|^2=\sum_{i=1}^d x_i^2$. 

We are interested in $\sup_{t\in[0,T]}{\cal W}_\rho({\cal L}(X_t),{\cal     L}(\bar{X}_t))$. Thanks to the Kantorovitch duality (see Corollary 2.5.2 in Rachev and R\"uschendorf~\cite{raru}), we know that for $t\in[0,T]$,
$${\cal W}_1({\cal L}(X_t),{\cal     L}(\bar{X}_t)) = \sup_{f:\R^d \rightarrow \R,\ Lip(f)\le 1} |\E[f(\bar{X}_t)-f(X_t)]|,$$
where $ Lip(f)=\sup_{x \not = y} \frac{|f(x)-f(y)|}{|x-y|}$. From the weak error expansion given by Talay and Tubaro~\cite{tt} when the coefficients are smooth enough, we deduce that ${\cal W}_1({\cal L}(X_T),{\cal     L}(\bar{X}_T))\ge  \frac{C}{N}$ for some constant $C>0$. Since, by H\"older's inequality, $\rho\mapsto{\cal W}_\rho$ is non-decreasing, we cannot therefore hope the order of convergence of $\sup_{t\in[0,T]}{\cal W}_\rho({\cal L}(X_t),{\cal     L}(\bar{X}_t))$ to be better than one. On the other hand, as remarked by Sbai~\cite{thesesbai}, a result of Gobet and Labart~\cite{goblab} supposing uniform ellipticity and some regularity on $\sigma$ and $b$ that will be made precise below implies that
$$ \sup_{t\in[0,T]}{\cal W}_1({\cal L}(X_t),{\cal     L}(\bar{X}_t))\leq \frac{C}{N}.$$

In a recent paper~\cite{AJK}, we proved that in dimension $d=1$, under uniform ellipticity and for coefficients $b$ and $\sigma$ time-homogeneous, bounded together with their derivatives up to the order $4$, one has \begin{equation}
   \sup_{t\in[0,T]}{\cal W}_\rho({\cal L}(X_t),{\cal     L}(\bar{X}_t))\leq \frac{C\sqrt{\ln (N)}}{N}\label{estiwas}
\end{equation} for any $\rho>1$. For the proof, we used that in dimension one, the optimal coupling measure $\pi$ between the measures $\mu$ and $\nu$ in the definition \eqref{defwas} of the Wasserstein distance is explicitly given by the inverse transform sampling: $\pi$ is the image of the Lebesgue measure on $[0,1]$ by the couple of pseudo-inverses of the cumulative distribution functions of $\mu$ and $\nu$. Our main result in the present paper is the generalization of \eqref{estiwas} to any dimension $d$ when the coefficients $b$ and $\sigma$ are time-homogeneous $C^2$, bounded together with their derivatives up to the order $2$ and uniform ellipticity holds.  We also generalize the analysis to time-dependent coefficients $b$ and $\sigma$ H\"older continuous with exponent $\gamma$ in the time variable. For $\gamma\in (0,1)$, the rate of convergence worsens i.e. the right-hand side of \eqref{estiwas} becomes $\frac{C}{N^\gamma}$ whereas it is preserved in the Lipschitz case $\gamma=1$.  These results are stated in Section \ref{sec:mainthm} together with the remark that the choice of a non-uniform time grid refined near the origin for the Euler scheme permits to get rid of the $\sqrt{\ln (N)}$ term in the numerator in the case $\gamma=1$. To our knowledge, they provide a new estimation of the weak error of the Euler scheme when the coefficients $b$ and $\sigma$ are only H\"older continuous in the time variable. The main difficulty to prove them is that, in contrast with the one-dimensional case, the optimal coupling between ${\cal L}(X_t)$ and ${\cal L}(\bar{X}_t)$ is only characterized in an abstract way. We want to apply the theory by Ambrosio {\it et al.} \cite{ags} to compute the time derivative $\frac{d}{dt}{\cal W}^\rho_\rho({\cal L}(X_t),{\cal     L}(\bar{X}_t))$. To do so, we have to interpret the Fokker-Planck equations giving the time derivatives of the densities of $X_t$ and $\bar{X}_t$ with respect to the Lebesgue measure as transport equations : the contribution of the Brownian term has to be written in the same way as the one of the drift term. This requires some regularity properties of the densities. In Section \ref{sec:prmainthm}, we give a heuristic proof of our main result without caring about these regularity properties. This allows us to present in a heuristic and pedagogical way the  main arguments, and to introduce the notations related to the optimal transport theory. In the obtained expression for $\frac{d}{dt}{\cal W}^\rho_\rho({\cal L}(X_t),{\cal     L}(\bar{X}_t))$, it turns out that, somehow because of the first order optimality condition on the optimal transport maps at time $t$, the derivatives of these maps with respect to the time variable do not appear (see Equation \eqref{der_Wass} below). The contribution of the drift term is similar to the one that we would obtain when computing $\frac{d}{dt}\E(|X_t-\bar{X}_t|^\rho)$ i.e. when working with the natural coupling between the SDE \eqref{sde} and its Euler scheme. To be able to deal with the contribution of the Brownian term, we first have to perform a spatial integration by parts. Then the uniform ellipticity condition enables us to apply a key lemma on pseudo-distances between matrices to see that this contribution is better behaved than the corresponding one in $\frac{d}{dt}\E(|X_t-\bar{X}_t|^\rho)$ and derive a stability inequality for ${\cal W}^\rho_\rho({\cal L}(X_t),{\cal     L}(\bar{X}_t))$ analogous to the one obtained in dimension $d=1$ in \cite{AJK}. Like in this paper, we conclude the heuristic proof by a Gronwall's type argument using estimations based on Malliavin calculus. In \cite{AJK}, our main motivation was to analyze the Wasserstein distance between the pathwise laws ${\cal L}((X_t)_{t\in [0,T]})$ and ${\cal L}((\bar{X}_t)_{t\in [0,T]})$. This gives then an upper bound of the error made when one approximates the expectation of a pathwise functional of the diffusion by the corresponding one computed with the Euler scheme. We were able to deduce from the upper bound on the Wasserstein distance between the marginal laws that the pathwise Wasserstein distance is upper bounded by~$CN^{-2/3+\varepsilon}$, for any $\varepsilon>0$. This improves the $N^{-1/2}$ rate given by the strong error analysis by Kanagawa~\cite{Ka}. To do so, we established using the Lamperti transform some key stability result for one-dimensional diffusion bridges in terms of the couple of initial and terminal positions. So far, we have not been able to generalize this stability result to higher dimensions. Nevertheless, our main result  can be seen as a first step in order to improve the estimation of the pathwise Wasserstein distance deduced from the strong error analysis.

In Section \ref{Sec_rig}, we give a rigorous proof of the main result. The theory of Ambrosio {\it et al.} \cite{ags} has been recently applied to Fokker-Planck equations associated with linear SDEs and SDEs nonlinear in the sense of McKean by Bolley {\it et al.}~\cite{bgg1,bgg2} in the particular case $\sigma=I_d$ of an additive noise and for the quadratic Wasserstein distance $\rho=2$ to study the long-time behavior of their solutions. In the present paper, we want to estimate the error introduced by a discretization scheme on a finite time-horizon with a general exponent $\rho$ and a non-constant diffusion matrix $\sigma$.  It turns out that, due to the local Gaussian behavior of the Euler scheme on each time-step, it is easier to apply the theory of Ambrosio {\it et al.} \cite{ags} to this scheme than to the limiting SDE \eqref{sde}. The justification of the spatial integration by parts performed on the Brownian contribution in the time derivative of the Wasserstein distance is also easier for the Euler scheme. That is why introduce a second Euler scheme
with time step $T/M$ and estimate the Wasserstein distance between the marginal laws of the two Euler schemes. We conclude the proof by letting $M\to\infty$ in this estimation thanks to the lower-semicontinuity of the Wasserstein distance with respect to the narrow convergence. The computation of the time derivative of the Wasserstein distance between the time-marginals of two Euler schemes can be seen as a first step to justify the formal expression of the time derivative of the Wasserstein distance between the time-marginals of the two limiting SDEs. We plan to investigate this problem in a future work. 

Section \ref{sectchlem} is devoted to technical lemmas including the already mentioned key lemma on the pseudo-distances between matrices and estimations based on Malliavin calculus.

\subsection*{Notations}
\begin{itemize}
   \item Unless explicitly stated, vectors are consider as column vectors.
\item The set of real $d\times d$ matrices is denoted by ${\cal M}_d(\R)$.
\item For a symmetric positive semidefinite matrix $M\in{\cal M}_d(\R)$, $M^{\frac{1}{2}}$ denotes the symmetric positive semidefinite matrix such that $M=M^{\frac{1}{2}}M^{\frac{1}{2}}$.
\item For $n\in \N$, we introduce
\begin{align*}
   C^{0,n}_b(\R) =\{& f:[0,T] \times \R^d\rightarrow \R  \text{ continuous, bounded and } \\&n \text{ times continuously differentiable in its $d$ last variables with bounded derivatives}\},
\end{align*}For $\gamma \in [0,1]$, we also define 
$$ C^{\gamma,n}_b(\R) =\{ f \in   C^{0,n}_b(\R), \text{ s. t. } \exists K\in[0,+\infty), \forall  s,t\in[0,T],\forall x \in \R^d, |f(t,x)-f(s,x)| \le K |t-s|^\gamma \},
$$ 
\begin{align*}
C^{\gamma,n}_b(\R^d) &=\{ f:[0,T]\times\R^d \rightarrow \R^d \text{ such that } \forall 1\le  i \le d, f_i\in C^{\gamma,n}_b(\R) \}, \\
C^{\gamma,n}_b(\mathcal{M}_d(\R)) &=\{ f:[0,T]\times\R^d \rightarrow \mathcal{M}_d(\R) \text{ such that } \forall 1\le  i,j \le d, f_{ij}\in C^{\gamma,n}_b(\R) \}.
\end{align*}
\item For $f:\R^d\to\R$ differentiable and $g:\R^d\to\R^d$, we denote by $\nabla f(g(x))$ the gradient $(\partial_{x_i}f)_{1\leq i\leq d}$ of $f$ computed at $g(x)$.
\item For $f:\R^d\to\R^d$, we denote by $\nabla f$ the Jacobian matrix $(\partial_{x_i}f_j)_{1\leq i,j\leq d}$ and by $\nabla^* f$ its transpose. 
\item For $f:\R^d\to\R$, we denote by $\nabla^2f$ the Hessian matrix $(\partial_{x_ix_j}f)_{1\leq i,j\leq d}$.
\item For $f:E\times \R^d\to\R$, we denote by $\nabla_x f(e,x)$, the partial gradient of $f$ with respect to its $d$ last variables.
\item For two density functions $p$ and $\bar{p}$ on $\R^d$, if there is a measurable function $f:\R^d\to\R^d$ such that the image of the probability measure $p(x)dx$ by $f$ admits the density $\bar{p}$, we write $p\#f=\bar{p}$.
\end{itemize}

\section{The main result}\label{sec:mainthm}
Our main result is the following theorem.
\begin{thm}\label{mainthm}
Assume that\begin{itemize}
\item $b \in C^{\gamma,2}_b(\R^d)$,
\item $\sigma \in C^{\gamma,2}_b(\mathcal{M}_d(\R))$ and is such that $a(t,x)=\sigma(t,x)\sigma(t,x)^*$ is uniformly elliptic, i.e.
$$\exists \underline{a}>0 \ s.t. \  \forall t\in[0,T],\;\forall x \in \R^d, \ a(t,x)-\underline{a} I_d \text{ is positive semidefinite}. $$
\end{itemize} 
Then
\begin{equation}
   \forall \rho \geq 1,\;\exists C<+\infty,\;\forall N\geq 1,\;\sup_{t\in[0,T]}{\cal W}_\rho({\cal L}(X_t),{\cal     L}(\bar{X}_t))\leq \frac{C \left(1+\mathbf{1}_{\gamma=1} \sqrt{\ln(N)}\right)}{N^\gamma},\label{resprinc}
\end{equation}
 where $C$ is a positive constant that only depends on $\rho$, $\underline{a}$,   $(\|\partial_{\alpha} a\|_\infty,\|\partial_{\alpha} b \|_\infty,0\le |\alpha|\le 2)$, and the coefficients $K,q$ involved in the $\gamma$-Hölder time regularity of $a$ and $b$. In particular $C$ does not depend on
  the initial condition~$x_0\in \R$.
\end{thm}
\begin{arem}\label{remssln}
Under the assumptions of Theorem \ref{mainthm} with $\gamma=1$, by discretizing the SDE \eqref{sde} with the Euler scheme on the non-uniform time grids refined near the origin $\left(t_i=(\frac{i}{N})^\beta T\right)_{0\leq i\leq N}$ with $\beta>1$, one gets rid of the $\sqrt{\ln(N)}$ term in the numerator : $$\exists C<+\infty,\;\forall N\geq 1,\;\sup_{t\in[0,T]}{\cal W}_\rho({\cal L}(X_t),{\cal     L}(\bar{X}_t))\leq \frac{C }{N}.$$
For $\gamma<1$, the choice of such non-uniform time grids does not lead to an improvement of the convergence rate in \eqref{resprinc}. For more details, see Remark \ref{justdisplnn} below.
\end{arem}
To our knowledge, Theorem \ref{mainthm} is a new result concerning the weak error of the Euler scheme, for coefficients $\sigma,b$ only $\gamma$-H\"older continuous in the time variable with $\gamma<1$. For $\gamma=1$, as remarked by Sbai~\cite{thesesbai}, a result of Gobet and Labart~\cite{goblab} supposing uniform ellipticity and that $b\in C^{1,3}_b(\R^d),\sigma \in  C^{1,3}_b(\mathcal{M}_d(\R))$ are continuously differentiable in time, implies that
$$ \sup_{t\in[0,T]}{\cal W}_1({\cal L}(X_t),{\cal     L}(\bar{X}_t))\leq \frac{C}{N}.$$Compared to this result, we have a slightly less accurate upper bound due to the $\sqrt{\ln(N)}$ term, but Theorem~\ref{mainthm} requires slightly less assumptions on the diffusion coefficients and most importantly concerns any $\rho$-Wasserstein distance. Using H\"older's inequality and the well-known boundedness of the moments of both $X_t$ and $\bar{X}_t$ for $t\in[0,T]$, one deduces that
\begin{acor}
   For any function $f:\R^d\to\R$ such that
$$\exists\alpha\in (0,1],\;\exists C,q\in(0,+\infty),\;\forall x,y\in\R^d,|f(x)-f(y)|\leq C(1+|x|^q+|y|^q)|x-y|^\alpha,$$
one has
\begin{align*}
 \exists C<+\infty,\;\forall N\geq 1,\;\sup_{t\in[0,T]}|\E(f(X_t))-\E(f(\bar{X}_t))|\leq \frac{C \left(1+\mathbf{1}_{\gamma=1} \sqrt{\ln(N)}\right)^\alpha}{N^{\alpha\gamma}}.
\end{align*}
\end{acor}
\begin{arem}\label{rem_polgrowth}
We have stated Theorem~\ref{mainthm} under assumptions that lead to a constant $C$  that does not depend on the initial condition~$x_0$. This is a nice feature that we used in~\cite{AJK} to bound the Wasserstein distance between the pathwise laws ${\cal L}((X_t)_{t\in [0,T]})$ and ${\cal L}((\bar{X}_t)_{t\in [0,T]})$ from above. However, Theorem~\ref{mainthm} still holds with a constant $C$ depending in addition on~$x_0$ if we relax the assumptions on $b$ and $\sigma$ as follows:
\begin{itemize}
\item $b$ and $\sigma$ are globally Lipschitz with respect to~$x$, i.e.
$$\forall f\in \{b,\sigma\}, \exists K\in[0,+\infty), \forall t\in[0,T], \forall  x,y \in \R^d, |f(t,y)-f(t,x)| \le K |x-y|,$$
\item $b$ and $\sigma$ are twice continuously differentiable in~$x$ and $\gamma$-H\"older in time, and such that we have the following polynomial growth
\begin{align*}
&\forall f\in \{b,\sigma\}, \exists K,q\in[0,+\infty),\; \forall  s,t\in[0,T], \forall x \in \R^d, |f(t,x)-f(s,x)| \le K |t-s|^\gamma(1+|x|^q), \\
&\text{for any } 1\le i,j,k,l\le d, \alpha \in \N^d,\text{ such that } |\alpha|=2 \text{ and } f\in \{\partial_{x_kx_l} b_i, \partial_{x_kx_l}\sigma_{ij} \}, \\
& \hspace{4cm}\exists K,q>0, \forall  t\ge 0, x \in \R^d, |f(t,x)| \le K (1+|x|^q),
\end{align*}
\item  $a(t,x)=\sigma(t,x)\sigma(t,x)^*$ is uniformly elliptic.
\end{itemize}
\end{arem}

Since by Hölder's inequality, $\rho\mapsto{\cal W}_\rho$ is non-increasing, it is sufficient to prove Theorem~\ref{mainthm} for $\rho$ large enough. In fact, we will assume through the rest of the article without loss of generality that $\rho\ge 2$. The main reason for this assumption is that the function $\R^d\times\R^d\ni(x,y)\mapsto|x-y|^\rho$, which appears in the definition \eqref{defwas} of $W_\rho$, becomes globally $C^2$. This will be convenient 
when studying the second order optimality condition. Furthermore, note that by the uniform ellipticity and regularity assumptions in Theorem \ref{mainthm}, for $t\in(0,T]$, $X_t$ and $\bar{X}_t$ admit densities respectively denoted by $p_t$ and $\bar{p}_t$ with respect to the Lebesgue measure. By a slight abuse of notation, we still denote by $W_\rho(p_t,\bar{p}_t)$ the $\rho$-Wasserstein distance between the probability measures $p_t(x)dx$ and $\bar{p}_t(x)dx$ on~$\R^d$.

\section{Heuristic proof of the main result}\label{sec:prmainthm}
The heuristic proof of Theorem \ref{mainthm} is structured as follows. First, we recall some optimal transport results about the Wasserstein distance and its associated optimal coupling, and we make some simplifying assumptions on the optimal transport maps that will be removed in the rigorous proof. Then, we can heuristically calculate $\frac{d}{dt}{\cal W}^\rho_\rho(p_t,\bar{p}_t)$, and get a sharp upper bound for this quantity. Last, we use a Gronwall's type argument to conclude the heuristic proof.

\subsection{Preliminaries on the optimal transport for the Wasserstein distance}\label{sec:characopttransp}

We introduce some notations that are rather standard in the theory of optimal transport (see \cite{ags,raru,villani}) and which will be useful to characterize the optimal coupling for the $\rho$-Wasserstein distance.  We will say that a function $\psi:\R^d\rightarrow [-\infty,+\infty]$ is $\rho$-convex if there is a function $\zeta:\R^d\rightarrow [-\infty,+\infty]$ such that 
$$\forall x \in \R^d, \  \psi(x)=\sup_{y\in \R^d}\left(-|x-y|^\rho-\zeta(y)\right).$$
In this case, we know from Proposition~3.3.5 of Rachev and R\"uschendorf~\cite{raru} that
\begin{equation}
   \forall x\in \R^d, \ \psi(x)=\sup_{y\in \R^d}\left( -|x-y|^\rho-\bar{\psi}(y)\right), \text{ where for } y\in \R^d, \bar{\psi}(y):=\sup_{x\in \R^d}\left( -|x-y|^\rho-\psi(x)\right).\label{rhoconv}
\end{equation}
We equivalently have,
\begin{equation}
  {\psi}(x)=-\inf_{y\in\R^d}\left(|x-y|^\rho+\bar{\psi}(y)\right) \mbox{ and } \bar{\psi}(x)=-\inf_{y\in\R^d}\left(|x-y|^\rho+\psi(y)\right).\label{dual}
\end{equation}
This result can be seen as an extension of the well-known Fenchel-Legendre duality for convex functions which corresponds to the case~$\rho=2$. We then introduce the $\rho$-subdifferentials of these functions. These are the sets defined by
\begin{align}
   \partial_\rho{\psi}(x)&=\{y\in\R^d:\psi(x)=-(|x-y|^\rho+\bar{\psi}(y))\},\label{diff1}\\
\partial_\rho\bar{\psi}(x)&=\{y\in\R^d:\bar{\psi}(x)=-(|x-y|^\rho+\psi(y)) \}.\label{diff2}
\end{align}

Let $t\in[0,T]$. According to Theorem 3.3.11 of Rachev and R\"uschendorf~\cite{raru}, we know that there is a couple~$(\xi_t,\bar{\xi}_t)$ of random variables with respective densities $p_t$ and $\bar{p}_t$ which attains the $\rho$-Wasserstein distance : $$ \E[|\xi_t-\bar{\xi}_t|^\rho]=W_\rho^\rho(p_t,\bar{p}_t).$$
Such a couple is called an optimal coupling for the Wasserstein distance. Besides, there exist two $\rho$-convex function~$\psi_t$ and $\bar{\psi}_t$ satisfying the duality property \eqref{rhoconv} and such that
$$\bar{\xi}_t \in \partial_\rho{\psi}_t( \xi_t) \text{ and } \xi_t \in \partial_\rho{\bar{\psi}}_t( \bar{\xi}_t), \ a.s.. $$

Now that we have recalled this well known result of optimal transport, we can start our heuristic proof of Theorem~\ref{mainthm}. To do so, we will assume that the  $\rho$-subdifferentials $\partial_\rho{\psi}_t(x)$ and $\partial_\rho\bar{\psi}_t(x)$ are non empty and single valued for any $x\in \R^d$, i.e. 
$$ \partial_\rho{\psi}_t(x)=\{ T_t(x) \}, \ \partial_\rho{\bar{\psi}}_t(x)=\{ \bar{T}_t(x) \}.$$
The functions $T_t(x)$ and $\bar{T}_t(x)$ depend on $\rho$ but we do not state explicitly this dependence for notational simplicity. Now, we clearly have
\begin{align}
   \psi_t(x)=-\left[|x-T_t(x)|^\rho+\bar{\psi}_t(T_t(x))\right] 
\mbox{and }\bar{\psi}_t(x)=-\left[|x-\bar{T}_t(x)|^\rho+{\psi}_t(\bar{T}_t(x))\right]. 
\label{exprexpl}
\end{align}
Besides, we can write the Wasserstein distance as follows:
\begin{equation}\label{Wass_expr}W_\rho^\rho(p_t,\bar{p}_t)=\int_{\R^d}|x-T_t(x)|^\rho p_t(x)dx=\int_{\R^d}|x-\bar{T}_t(x)|^\rho \bar{p}_t(x)dx.
\end{equation}
Since on the one hand $\bar{\xi}_t = T_t(\xi_t)$ and $\xi_t=\bar{T}_t(\bar{\xi}_t)$ almost surely, and on the other hand $p_t(x)\bar{p}_t(x)>0$ thanks to the uniform ellipticity assumption,
 \begin{equation}
dx\mbox{ a.e.},\;\bar{T}_t( T_t(x) )=T_t(\bar{T}_t(x))=x.\label{inv}
 \end{equation}
In the remaining of Section \ref{sec:prmainthm}, we will perform heuristic computations without caring about the actual smoothness of the functions ${\psi}_t$, ${\bar{\psi}}_t$, $T_t$ and $\bar{T}_t$.  In particular, we suppose that
\begin{align}
\forall x \in \R^d, \bar{T}_t( T_t(x) )&=T_t(\bar{T}_t(x))=x \label{invpsi}\\\nabla\bar{\psi}_t(T_t(x))&=\rho|x-T_t(x)|^{\rho-2}(x-T_t(x)),\label{eul1}\\
\nabla {\psi}_t(\bar{T}_t(x))&=\rho|x-\bar{T}_t(x)|^{\rho-2}(x-\bar{T}_t(x)).\label{eul2}
\end{align}
where the two last equations are the first order Euler conditions of optimality in the minimization problems \eqref{dual}. 

\subsection{A formal computation of  $\frac{d}{dt}W_\rho^\rho(p_t,\bar{p}_t)$}\label{sec:formderw}
We now make a heuristic differentiation of~\eqref{Wass_expr} with respect to~$t$. A computation of the same kind for the case $\rho=2$ and with identity diffusion matrix $\sigma$ is given by Bolley {\it et al.} : see  p.2437 and Remark 3.6 p.2445 in~\cite{bgg1} or p.431 in~\cite{bgg2}. 
\begin{align*}
  \frac{d}{dt}W_\rho^\rho(p_t,\bar{p}_t) & =\int_{\R^d}|x-T_t(x)|^\rho\partial_tp_t(x)dx+\int_{\R^d}\rho|x-T_t(x)|^{\rho-2}(T_t(x)-x).\partial_tT_t(x)p_t(x)dx\\
&=\int_{\R^d}|x-T_t(x)|^\rho\partial_tp_t(x)dx-\int_{\R^d}\nabla\bar{\psi}_t(T_t(x)).\partial_tT_t(x)p_t(x)dx\\
&=\int_{\R^d}\left(|x-T_t(x)|^\rho+\bar{\psi}_t(T_t(x))\right)\partial_tp_t(x)dx \\&-\int_{\R^d}\left(\nabla\bar{\psi}_t(T_t(x)).\partial_tT_t(x)p_t(x)+\bar{\psi}_t(T_t(x))\partial_tp_t(x)\right)dx\\
&=-\int_{\R^d}\psi_t(x)\partial_tp_t(x)dx-\frac{d}{dt}\int_{\R^d}g(T_t(x))p_t(x)dx\bigg|_{g=\bar{\psi}_t},
\end{align*}
where we used~\eqref{eul1} for the second equality and~\eqref{exprexpl} for the fourth.
Since the image of the probability measure $p_t(x)dx$ by the map $T_t$ is the probability measure $\bar{p}_t(x)dx$, which we write as $$\bar{p}_t=T_t\#p_t$$ in what follows, we have $\int_{\R^d}g(T_t(x))p_t(x)dx=\int_{\R^d}g(x)\bar{p}_t(x)dx$ and thus $\frac{d}{dt}\int_{\R^d}g(T_t(x))p_t(x)dx=\int_{\R^d}g(x) \partial_t \bar{p}_t(x)dx$. This heuristic calculation finally gives
\begin{equation}\label{der_Wass}
\frac{d}{dt}W_\rho^\rho(p_t,\bar{p}_t)=-\int_{\R^d}\psi_t(x)\partial_tp_t(x)dx-\int_{\R^d}\bar{\psi}_t(x)\partial_t\bar{p}_t(x)dx.
\end{equation}

Let us assume now that the following Fokker-Planck equations  for the densities $p_t$ and $\bar{p}_t$ hold in the classical sense
\begin{align}
 &\partial_t p_t(x)=\frac{1}{2}\sum_{i,j=1}^d\partial_{x_ix_j}(a_{ij}(t,x)p_t(x))-\sum_{i=1}^d\partial_{x_i}(b_i(t,x)p_t(x)),\label{fpp}\\
&\partial_t \bar{p}_t(x)=\frac{1}{2}\sum_{i,j=1}^d\partial_{x_ix_j}(\bar{a}_{ij}(t,x)\bar{p}_t(x))-\sum_{i=1}^d\partial_{x_i}(\bar{b}_i(t,x)\bar{p}_t(x)),
\label{fpbarp}\end{align}
where 
\begin{equation}\label{def_aetb_bar}\left(\begin{array}{c}\bar{a}
   \\\bar{b}
\end{array}\right)(t,x)=\E\left(\left(\begin{array}{c}a
   \\b
\end{array}\right)(\tau_t,\bar{X}_{\tau_t})|\bar{X}_t=x\right).
\end{equation}
The first equation is the usual Fokker-Planck equation for the SDE \eqref{sde}. For the second one, we also use the result by Gy\"ongy~\cite{gyongy} that ensures that the SDE with coefficients $\bar{b}$ and $\bar{a}^{\frac{1}{2}}$ has the same marginal laws as the Euler scheme. Now, plugging these equations in~\eqref{der_Wass}, we get 
\begin{align}
   \frac{d}{dt}W_\rho^\rho(p_t,\bar{p}_t)=&-\frac{1}{2}\int_{\R^d}\Tr(\nabla^2\psi_t(x)a(t,x))p_t(x)dx-\int_{\R^d}\nabla\psi_t(x).b(t,x)p_t(x)dx\notag\\&-\frac{1}{2}\int_{\R^d}\Tr(\nabla^2\bar{\psi}_t(x)\bar{a}(t,x))\bar{p}_t(x)dx-\int_{\R^d}\nabla\bar{\psi}_t(x).\bar{b}(t,x)\bar{p}_t(x)dx\notag
\end{align}
by using integrations by parts and assuming that the boundary terms vanish. We now use $\bar{p}_t = T_t\#p_t$ and 
\begin{equation}
   \nabla{\psi}_t(x)=\rho|T_t(x)-x|^{\rho-2}(T_t(x)-x)=-\nabla\bar{\psi}_t(T_t(x)),\label{eul1comp}
\end{equation}
which is deduced from~\eqref{invpsi},~\eqref{eul1}  and~\eqref{eul2}, to get
\begin{align}
 \frac{d}{dt}W_\rho^\rho(p_t,\bar{p}_t)=&-\frac{1}{2}\int_{\R^d}\Tr[\nabla^2\psi_t(x)a(t,x)+\nabla^2\bar{\psi}_t(T_t(x)) \bar{a}(t,T_t(x)) ]p_t(x)dx\notag\\
&-\int_{\R^d}\left(\nabla\psi_t(x).b(t,x)+\nabla\bar{\psi}_t(T_t(x)).\bar{b}(t,T_t(x))\right)p_t(x)dx\notag\\
=&-\frac{1}{2}\int_{\R^d}\Tr[\nabla^2\psi_t(x)a(t,x)+\nabla^2\bar{\psi}_t(T_t(x))\bar{a}(t,T_t(x))]p_t(x)dx\notag\\&+\rho\int_{\R^d}|T_t(x)-x|^{\rho-2}(T_t(x)-x).\left(\bar{b}(t,T_t(x))-b(t,x)\right)p_t(x)dx\label{evolwass}.
\end{align}
This formula looks very nice but due to the lack of regularity of $\psi_t$ and $\bar{\psi}_t$, which are merely semi-convex functions, it is only likely to hold with the equality replaced by $\leq$ and the $\nabla^2\psi_t$ and $\nabla^2\bar{\psi}_t$ replaced by the respective Hessians in the sense of Alexandrov of $\psi_t$ and $\bar{\psi}_t$. See Proposition \ref{prop_1} where such an inequality is proved rigorously for the Wassertein distance between the time marginals of two Euler schemes.

\subsection{Derivation of a stability inequality for $W_\rho^\rho(p_t,\bar{p}_t)$}\label{sec:stabineq}
In \eqref{evolwass}, the contribution of the drift terms only involves the optimal transport and is equal to $\rho\E\left(|\bar{\xi}_t-\xi_t|^{\rho-2}(\bar{\xi}_t-\xi_t).\left(\bar{b}(t,\bar{\xi}_t)-b(t,\xi_t)\right)\right)$ for any optimal coupling $(\xi_t,\bar{\xi}_t)$ between $p_t$ and $\bar{p}_t$. To obtain this term, it was enough to use the first order optimality conditions \eqref{eul1} and \eqref{eul2}. To deal with the Hessians $\nabla^2\psi_t$ and $\nabla^2\bar{\psi}_t$ which appear in the contribution of the diffusion terms, we will need the associated second order optimality conditions.

Differentiating~\eqref{eul1comp} with respect to~$x$, we get
\begin{equation}
\nabla^2\psi_t(x)=\rho|T_t(x)-x|^{\rho-2}\left(I_d+(\rho-2)\frac{T_t(x)-x}{|T_t(x)-x|}\frac{(T_t(x)-x)^*}{|T_t(x)-x|}\right)\left(\nabla^*T_t(x)- I_d\right).
\label{hesspsi}\end{equation}
By symmetry and \eqref{invpsi},
$$\nabla^2\bar{\psi}_t(T_t(x))=\rho|T_t(x)-x|^{\rho-2}\left(I_d+(\rho-2)\frac{T_t(x)-x}{|T_t(x)-x|}\frac{(T_t(x)-x)^*}{|T_t(x)-x|}\right)\left(\nabla^*\bar{T}_t(T_t(x))-I_d\right).$$
By differentiation of \eqref{invpsi}, we get that $\nabla^*T_t(x)$ is invertible and have $\nabla^*\bar{T}_t(T_t(x))=(\nabla^*T_t(x))^{-1}$. Plugging these equations into~\eqref{evolwass}, we get 
\begin{align*}
\frac{d}{dt}W_\rho^\rho(p_t,\bar{p}_t)   =&\rho\int_{\R^d}|T_t(x)-x|^{\rho-2}(T_t(x)-x).\left(\bar{b}(t,T_t(x))-b(t,x)\right)p_t(x)dx \\ &+\frac{\rho}{2}\int_{\R^d}|T_t(x)-x|^{\rho-2} \Tr\bigg[\bigg(I_d+(\rho-2)\frac{T_t(x)-x}{|T_t(x)-x|}\frac{(T_t(x)-x)^*}{|T_t(x)-x|}\bigg)\\ &\left\{\left(I_d-\nabla^*T_t(x)\right)a(t,x)+\left(I_d-(\nabla^*T_t(x))^{-1}\right)\bar{a}(t,T_t(x))\right\}\bigg]p_t(x)dx.
\end{align*}
In order to make the diffusion contribution of the same order as the drift one, we want to upper-bound the trace term by the square of a distance between $a(t,x)$ and $\bar{a}(t,T_t(x))$. The key Lemma~\ref{lemma_Wp2} permits to do so. To check that its hypotheses are satisfied, we remark that the second order optimality condition for~\eqref{dual} 
\begin{equation}
\nabla^2\psi_t(\bar T_t(y))+\rho|\bar T_t(y)-y|^{\rho-2}\bigg(I_d+(\rho-2)\frac{\bar T_t(y)-y}{|\bar T_t(y)-y|}\frac{(\bar T_t(y)-y)^*}{|\bar T_t(y)-y|}\bigg)\geq 0\label{secondordeul}
\end{equation}
computed at $y=T_t(x)$ combined with \eqref{invpsi} and \eqref{hesspsi} gives that $$M:=\bigg(I_d+(\rho-2)\frac{T_t(x)-x}{|T_t(x)-x|}\frac{(T_t(x)-x)^*}{|T_t(x)-x|}\bigg)\nabla^*T_t(x)$$ is a positive semidefinite matrix. It is in fact positive since it is the product of two invertible matrices. We can then apply the key Lemma~\ref{lemma_Wp2} with $v=\frac{T_t(x)-x}{|T_t(x)-x|}$, $a_1=a(t,x)$, $a_2=\bar{a}(t,x)$ and $M$ defined just above and get:
\begin{align*}
\frac{d}{dt}W_\rho^\rho(p_t,\bar{p}_t)   \le&\rho\int_{\R^d}|T_t(x)-x|^{\rho-1}\left|b(t,x)-\bar{b}(t,T_t(x))\right|p_t(x)dx \\ &+ \frac{ \rho(\rho-1)^2 }{8 \underline{a} }  \int_{\R^d}|T_t(x)-x|^{\rho-2} \Tr\left[\left(a(t,x)-\bar{a}(t,T_t(x)) \right)^2 \right] p_t(x)dx. \nonumber
\end{align*}
Finally, using that $\bar{p}_t=T_t\#p_t$, we get
\begin{align}
\frac{d}{dt}W_\rho^\rho(p_t,\bar{p}_t)   \le&\rho\int_{\R^d}|x-\bar{T}_t(x)|^{\rho-1}\left|b(t,\bar{T}_t(x))-\bar{b}(t,x)\right| \bar{p}_t(x)dx \label{maj_evolwass}\\ &+\frac{ \rho(\rho-1)^2 }{8 \underline{a} }   \int_{\R^d}|x-\bar{T}_t(x)|^{\rho-2} \Tr\left[\left(a(t,\bar{T}_t(x))-\bar{a}(t,x) \right)^2 \right]\bar{p}_t(x)dx. \nonumber
\end{align}
Now, we use the triangle inequalities $|b(t,\bar{T}_t(x))-\bar{b}(t,x)|\le |b(t,\bar{T}_t(x))-b(t,x)|+|b(t,x)-\bar{b}(t,x)|$ and 
$\Tr\left[\left(a(t,\bar{T}_t(x))-\bar{a}(t,x) \right)^2 \right]\le 2 \left[\Tr\left[\left(a(t,\bar{T}_t(x))-a(t,x) \right)^2 \right]+ \Tr\left[\left(a(t,x)-\bar{a}(t,x) \right)^2 \right] \right]$ together with the assumptions on $a$ and $b$ to get that there is a constant $C$ depending only on $\rho$, $\underline{a}$ and the spatial Lipschitz constants of~$a$ and $b$ such that
\begin{align}
\frac{d}{dt}W_\rho^\rho(p_t,\bar{p}_t)   \le& C \Bigg( W_\rho^\rho(p_t,\bar{p}_t)+\int_{\R^d}|x-\bar{T}_t(x)|^{\rho-1}\left|b(t,x)-\bar{b}(t,x)\right| \bar{p}_t(x)dx \label{maj_evolwass2}\\ &+  \int_{\R^d}|x-\bar{T}_t(x)|^{\rho-2} \Tr\left[\left(a(t,x)-\bar{a}(t,x) \right)^2 \right]\bar{p}_t(x)dx \Bigg). \nonumber
\end{align}

\begin{arem}
Equation~\eqref{maj_evolwass} illustrates the difference between the weak error and the strong error analysis. To study the strong error between $X_t$ and $\bar{X}_t$, one would typically apply Itô's formula and take expectations to get
\begin{align}
 \frac{d}{dt}&\E(|X_t-\bar{X}_t|^\rho)=\rho\E\Bigg(|X_t-\bar{X}_t|^{\rho-2}(X_t-\bar{X}_t).(b(t,X_t)-b(\tau_t,\bar{X}_{\tau_t}))+ \frac{1}{2}|X_t-\bar{X}_t|^{\rho-2}\notag\\&\times\Tr\bigg[\bigg(I_d+(\rho-2)\frac{X_t-\bar{X}_t}{|X_t-\bar{X}_t|}\frac{(X_t-\bar{X}_t)^*}{|X_t-\bar{X}_t|}\bigg)\times(\sigma(t,X_t)-\sigma(\tau_t,\bar{X}_{\tau_t}))(\sigma(t,X_t)-\sigma(\tau_t,\bar{X}_{\tau_t}))^*\bigg]\Bigg)\notag\\
&\leq C\E\Bigg(|X_t-\bar{X}_t|^{\rho-2}\bigg\{|X_t-\bar{X}_t|^{2}+(X_t-\bar{X}_t).(b(t,\bar{X}_t)-b(\tau_t,\bar{X}_{\tau_t}))\notag\\&+\Tr\bigg[\left(I_d+(\rho-2)\frac{X_t-\bar{X}_t}{|X_t-\bar{X}_t|}\frac{(X_t-\bar{X}_t)^*}{|X_t-\bar{X}_t|}\right)\times(\sigma(t,\bar{X}_t)-\sigma(\tau_t,\bar{X}_{\tau_t}))(\sigma(t,\bar{X}_t)-\sigma(\tau_t,\bar{X}_{\tau_t}))^*\bigg]\bigg\}\Bigg).\label{errfort}
\end{align}
The diffusion contribution is very different from the one in~\eqref{maj_evolwass2} : indeed, the absence of conditional expectation in the quadratic factor $(\sigma(t,\bar{X}_t)-\sigma(\tau_t,\bar{X}_{\tau_t}))(\sigma(t,\bar{X}_t)-\sigma(\tau_t,\bar{X}_{\tau_t}))^*$ in the trace term does not permit cancellations like in~\eqref{maj_evolwass2} where $\int_{\R^d}\Tr\left[\left(a(t,x)-\bar{a}(t,x) \right)^2 \right]\bar{p}_t(x)dx=\E\left(\Tr((\E(a(t,\bar{X}_t)-a(\tau_t,\bar{X}_{\tau_t})|\bar{X}_t))^2)\right)$.  

As an aside remark, we see that when $\sigma$ is constant, the diffusion contribution disappears in Equation \eqref{errfort} and is non-positive in Equation \eqref{evolwass}. In this case, $\sup_{t\in [0,T]}\E^{1/\rho}(|X_t-\bar{X}_t|^\rho)$ can be upper bounded by $C/N^\gamma$ where $\gamma$ denotes the H\"older exponent of the coefficient $b$ in the time variable. For $\gamma=1$, this leads to the improved bound $\sup_{t\in [0,T]}W_\rho(p_t,\bar{p}_t)\leq C/N$.
\end{arem}

\subsection{The argument based on Gronwall's lemma}\label{subs_Gronwall}

Starting from \eqref{maj_evolwass2}, we can conclude by applying a rigorous Gronwall type argument, which is analogous to the one used in the one-dimensional case in~\cite{AJK}. For the sake of completeness, we nevertheless repeat these calculations since we consider here in addition  coefficients which are not time-homogeneous but $\gamma$-Hölder continuous in time.

We set $\zeta_\rho(t)={\cal W}^2_\rho(p_t,\bar{p}_t)$ and define for any integer $k\geq 1$,
\begin{equation*}
   \;h_k(x)=k^{-2/\rho}h(kx)\mbox{ where }h(x)=\begin{cases}x^{2/\rho}\mbox{ if }x\geq 1,\\
1+\frac{2}{\rho}(x-1)\mbox{ otherwise.}
\end{cases}
\end{equation*}
Since $h_k$ is $C^1$ and non-decreasing, we get from~\eqref{maj_evolwass2} and H\"older's inequality 
\begin{align*}
   h_k\left(\zeta^{\rho/2}_\rho(t)\right)
&\leq h_k(0)+C\int_0^th_k'\left(\zeta^{\rho/2}_\rho(s)\right)\bigg[\zeta^{\rho/2}_\rho(s) \\ 
&\phantom{h_k(0)+C\int_0^t} +\zeta^{(\rho-1)/2}_\rho(s)\left( \int_{\R^d} |b(s,x)-\bar{b}(s,x)|^\rho \bar{p}_s(x)dx\right)^{1/\rho} \\
&\phantom{h_k(0)+C\int_0^t}  +\zeta^{(\rho-2)/2}_\rho(s)\left(\int_{\R^d} \Tr\left[\left(a(s,x)-\bar{a}(s,x) \right)^2 \right]^{\rho/2} \bar{p}_s(x)dx\right)^{2/\rho}\bigg]ds.
\end{align*}
Since $(h'_k)_{k\geq 1}$ is a non-decreasing sequence of functions that converges to $x\mapsto \frac{2}{\rho}x^{\frac{2}{\rho}-1}$ as $k\to\infty$, we get by the monotone convergence theorem and~\eqref{def_aetb_bar}
\begin{align*}\zeta_\rho(t)\le & \frac{2C}{\rho} \int_0^t \zeta_\rho(s)+ \zeta_\rho(s)^{1/2} \E^{1/\rho}[|b(s,\bar{X}_s)-\E[b(\tau_s,\bar{X}_{\tau_s})|\bar{X}_s ]|^\rho]  \\
&+\E^{2/\rho}[\Tr((a(s,\bar{X}_s)-\E[a(\tau_s,\bar{X}_{\tau_s})|\bar{X}_s ])^2)^{\rho/2}] ds. 
\end{align*}
Let us focus for example on the diffusion term. First,
$$\Tr[(a(s,\bar{X}_s)-\E[a(\tau_s,\bar{X}_{\tau_s})|\bar{X}_s ])^2]^{\rho/2}\le d^{\rho-2}\sum_{i,j=1}^d|a_{ij}(s,\bar{X}_s)-\E[a_{ij}(\tau_s,\bar{X}_{\tau_s})|\bar{X}_s ]|^\rho.$$ We have 
$$|a_{ij}(s,\bar{X}_s)-a_{ij}(\tau_s,\bar{X}_s)| \le K |s-\tau_s|^\gamma $$
and 
\begin{align}
a_{ij}(\tau_s,\bar{X}_s)-a_{ij}(\tau_s,\bar{X}_{\tau_s}) =&(\bar{X}_{s}-\bar{X}_{\tau_s}). \int_0^1 \nabla_x a_{ij}(\tau_s,v\bar{X}_{s}+(1-v)\bar{X}_{\tau_s} )dv \notag \\
=& \nabla_x a_{ij}(\tau_s,\bar{X}_{s}).\left[\sigma(\tau_s,\bar{X}_{s})(W_{s}-W_{\tau_s})\right] \notag\\ 
&+ \nabla_x a_{ij}(\tau_s,\bar{X}_{s}).\left[ (\sigma(\tau_s,\bar{X}_{\tau_s})-\sigma(\tau_s,\bar{X}_{s}))(W_s-W_{\tau_s}) + b(\tau_s,\bar{X}_{\tau_s})(s-\tau_s)\right] \notag \\
&+(\bar{X}_{s}-\bar{X}_{\tau_s}). \int_0^1 \nabla_x a_{ij}(\tau_s,v\bar{X}_{s}+(1-v)\bar{X}_{\tau_s} )- \nabla_x a_{ij}(\tau_s,\bar{X}_{s} )dv.\label{majodiffa}
\end{align}
Now, we use Jensen's inequality together with the boundedness of $b$ and the boundedness and Lipschitz property of~$x \mapsto \nabla_x a_{ij}(t,x)$, uniformly in $t\in[0,T]$, to get 
\begin{align*}
\E\left[|a_{ij}(s,\bar{X}_s)-\E[a_{ij}(\tau_s,\bar{X}_{\tau_s})|\bar{X}_s ]|^\rho \right] &\le  
\frac{C}{N^{\gamma \rho}}+ C \E[|\sigma^*(\tau_s,\bar{X}_s)\nabla_x a_{ij}(\tau_s,\bar{X}_{s})|^\rho|\E[(W_{\tau_s}-W_s)|\bar{X}_s] |^\rho ] \\
&  +C\left[\frac{1}{N^{\rho}}+\E[|(\sigma(\tau_s,\bar{X}_{\tau_s})-\sigma(\tau_s,\bar{X}_{s}))(W_s-W_{\tau_s})|^\rho+|\bar{X}_{\tau_s}-\bar{X}_{s}|^{2\rho}]\right].
\end{align*}
By the boundedness of $\sigma$ and $b$, one easily checks that
\begin{equation}
   \forall q\geq 1,\;\exists C\in[0,+\infty),\;\forall 0\leq s\leq t\leq T,\;\E(|\bar{X}_t-\bar{X}_s|^q)\leq C(t-s)^{q/2}.\label{majaccroiss}
\end{equation}
With Lemma~\ref{malcal} and the spatial Lipschitz continuity of $\sigma$, we deduce that
\begin{align}
\E\left[|a_{ij}(s,\bar{X}_s)-\E[a_{ij}(\tau_s,\bar{X}_{\tau_s})|\bar{X}_s ]|^\rho \right] &\le C(s-\tau_s)^{\gamma\rho}+C\left((s-\tau_s)\wedge \left (\frac{(s-\tau_s)^2}s+\frac{1}{N^2}\right )\right)^{\rho/2}\label{majofine}\\
&\leq \frac{C}{N^{\gamma \rho}}+\frac{C}{N^{\rho/2}\vee (N^\rho s^{\rho/2})}.\notag\end{align}
As a similar bound holds for the drift contribution, we finally get:
\begin{align*}\zeta_\rho(t)\le & C \int_0^t \zeta_\rho(s)+ \zeta_\rho(s)^{1/2} \left(\frac{1}{N^{\gamma }}+\frac{1}{N^{1/2} \vee (N s^{1/2}) } \right) +\frac{1}{N^{2 \gamma }}+\frac{1}{N \vee (N^2 s) } ds \\
&\le C \int_0^t \zeta_\rho(s) +\frac{1}{N^{2 \gamma }}+\frac{1}{N \vee (N^2 s) } ds \\
&\le C \int_0^t \zeta_\rho(s)ds + C\left( \frac{1}{N^{2\gamma }} + \frac{\ln(N)}{N^2}\right),
\end{align*}
and we obtain Theorem~\ref{mainthm} by Gronwall's lemma. 
\begin{arem}\label{justdisplnn}
In case $\gamma=1$, choosing $\beta>1$ and replacing the uniform time-grid by the grid $\left(t_i=(\frac{i}{N})^\beta T\right)_{0\leq i\leq N}$ refined near the origin, one may take advantage of \eqref{majofine} which is still valid with the last discretization time $\tau_t$ before $t$ now equal to $\left(\frac{\lfloor N (t/T)^{1/\beta}\rfloor}{N}\right)^\beta T$, since the largest step in the grid is $t_N-t_{N-1}\leq \frac{\beta T}{N}$. Adapting the above argument based on Gronwall's lemma, one obtains the statement in Remark \ref{remssln}. Indeed,  one has
\begin{align*}
   \int_0^T&\left(\frac{(s-\tau_s)^2}{s}+\frac{1}{N^2}\right)\wedge (s-\tau_s)ds
\leq \int_0^{T/N^\beta}(s-\tau_s)ds+\int_{T/N^\beta}^T\frac{(s-\tau_s)^2}{s}ds+\frac{T}{N^2}\\&=\frac{T^2}{2N^{2\beta}}+T^2\sum_{k=1}^{N-1}\left(\frac{k}{N}\right)^{2\beta}\left[\frac{1}{2}(1+1/k)^{2\beta}-\frac{1}{2}+2-2(1+1/k)^\beta+\beta\ln(1+1/k)\right]+\frac{T}{N^2}.
\end{align*}
Expanding the term between square brackets in powers of $1/k$, one easily checks that this term behaves like ${\mathcal O}(k^{-3})$. Now $$\sum_{k=1}^{N-1}\left(\frac{k}{N}\right)^{2\beta}\frac{1}{k^3}=N^{-2\beta}\sum_{k=1}^{N-1}k^{2\beta-3}=N^{-2\beta}{\mathcal O}(N^{2\beta-2})={\mathcal O}(N^{-2}).$$
One concludes that $$\exists C<+\infty,\;\forall N\geq 1,\;\int_0^T\left(\frac{(s-\tau_s)^2}{s}+\frac{1}{N^2}\right)\wedge (s-\tau_s)ds\leq \frac{C}{N^2}.$$
\end{arem}
\begin{arem}
If we only use the assumptions of Remark~\ref{rem_polgrowth}, we now deduce from \eqref{majodiffa} the existence of finite constants $C,q>0$ depending on~$\rho$,
\begin{align*}
\E\big[|a_{ij}(s,\bar{X}_s)&-\E[a_{ij}(\tau_s,\bar{X}_{\tau_s})|\bar{X}_s ]|^\rho \big] \le
\frac{C \E[(1+|\bar{X}_s|^q)]}{N^{\gamma \rho}}+ C \E[|\E[(W_{\tau_s}-W_s)|\bar{X}_s] |^\rho (1+|\bar{X}_s|^q)] \\
&  +C \E\Bigg[ (1+|\bar{X}_s|^q+|\bar{X}_{\tau_s}|^q) \left(\frac{1}{N^{\rho}} +|\bar{X}_{\tau_s}-\bar{X}_{s}|^\rho|W_s-W_{\tau_s}|^\rho + |\bar{X}_{\tau_s}-\bar{X}_{s}|^{2\rho} \right)  \Bigg].
\end{align*}
We can conclude that \eqref{resprinc} still holds with a constant $C$ depending on $x_0$ by using that the moments of the Euler scheme are uniformly bounded i.e. $\forall q'\geq 1$, $\E[\sup_{t\in [0,T]} |\bar{X}_t|^{q'} ]\le K_{q'}(1+|x_0|^{q'})$, an adaptation of Lemma~\ref{malcal} and the Cauchy-Schwarz inequality.
\end{arem}

\section{A rigorous proof of Theorem~\ref{mainthm}}\label{Sec_rig}

\subsection{Discussion of the strategy of the proof}\label{Subsec_rig}

We start by listing the simplifying hypotheses that we made in the Sections \ref{sec:characopttransp}, \ref{sec:formderw} and \ref{sec:stabineq}.
\begin{enumerate}
\item The $\rho$-subdifferentials $\partial_\rho \psi_t(x)$ and $\partial_\rho \bar{\psi}_t(x)$ are single valued.
\item The optimal transport and the densities~$p_t$ and $\bar{p}_t$ are smooth enough to get the time derivative of the Wasserstein distance~\eqref{der_Wass}.
\item The Fokker-Planck equations~\eqref{fpp} and \eqref{fpbarp} hold in the classical sense.
\item The functions $\psi_t$ and $\bar{\psi}_t$ are smooth enough and the integration by parts leading to~\eqref{evolwass} are valid.
\end{enumerate}
Let us now comment how we will manage to prove our main result without using these simplifying hypotheses. The first one was mainly used to get that the optimal transport maps are inverse functions (see~\eqref{invpsi} above). Still, the optimal transport theory will give us the existence of optimal transport maps that are inverse functions of each other. 

The second point is more crucial and is related to the third. Let us assume that there are Borel vector fields $v_t(x)$  and $\bar{v}_t(x)$ such that
\begin{equation}
   \int_0^T \left( \int_{\R^d} |v_t(x)|^\rho p_t(x)dx\right)^{1/\rho}dt+\int_0^T \left( \int_{\R^d} |\bar{v}_t(x)|^\rho \bar{p}_t(x)dx\right)^{1/\rho}dt<\infty\label{integv}
\end{equation}
and the so-called transport equations
\begin{equation}\partial_t p_t + \nabla .(v_t p_t)=0\mbox{ and }\partial_t \bar{p}_t + \nabla .(\bar{v}_t \bar{p}_t)=0\label{fp2}
\end{equation}
hold in the sense of distributions.
This means that for any $C^\infty$ function $\varphi$ with compact support on $(0,T) \times \R^d$,
$$\int_0^T\int_{\R^d} \left(\partial_t \varphi(t,x) + v_t(x).\nabla \varphi(t,x)  \right)p_t(x)dx dt=0,$$
and the same for $\bar{p}_t$. Then, it can be deduced from Ambrosio, Gigli and Savaré~\cite{ags} that $t\mapsto W_\rho^\rho(p_t,\bar{p}_t)$ is absolutely continuous and such that $dt$ a.e.,
\begin{equation}\label{der_ags1}
\frac{d}{dt} W_\rho^\rho(p_t,\bar{p}_t)=\rho\int_{\R^d}|T_t(x)-x|^{\rho-2}(x-T_t(x)).v_t(x)p_t(x)+|\bar{T}_t(x)-x|^{\rho-2}(x-\bar{T}_t(x)).\bar{v}_t(x)\bar{p}_t(x)dx.
\end{equation}
For the details, see the second paragraph called ``The time derivative of the Wasserstein distance'' in Subsection~\ref{subsec_prop44}.

Thus, it would be sufficient to show that the Fokker-Planck equations may be reformulated as the transport equations \eqref{fp2}. Concerning $p_t$, for the integrability condition \eqref{integv} to be satisfied by the natural choice $v_t(x)=b(t,x)-\frac{\nabla_x^*.(a(t,x)p_t(x))}{2p_t(x)}$ deduced from \eqref{fpp}, one typically needs \begin{equation}\label{integderlogrho}\int_0^T\left(\int_{\R^d} |\nabla_x \ln p_t(x)|^\rho p_t(x)dx\right)^{1/\rho}dt<+\infty
\end{equation}
For $\rho=2$, one may generalize the argument given by Bolley {\it et al.} p.2438~\cite{bgg1} in the particular case $\sigma=I_d$. Using \eqref{fpp} and an integration by parts for the last equality, one obtains formally
\begin{align*}
   \frac{d}{dt}&\int_{\R^d}\ln p_t(x)p_t(x)dx=\int_{\R^d}\ln p_t(x)\partial_tp_t(x)dx+\int_{\R^d}\partial_tp_t(x)dx\\&=\int_{\R^d}b(t,x).\nabla_xp_t(x)-\frac{1}{2}\sum_{i,j=1}^d\left(\partial_{x_i}p_t(x)\partial_{x_j}a_{ij}(t,x)+a_{ij}(t,x)\frac{\partial_{x_i}p_t(x)\partial_{x_j}p_t(x)}{p_t(x)}\right)dx+0
\end{align*}
to deduce with the uniform ellipticity condition and the positivity of the relative entropy $\int_{\R^d}\ln((2\pi)^{d/2}p_T(x)e^{|x|^2/2})p_T(x)dx$ that for $t_0\in (0,T]$,
\begin{align*}
   \int_{t_0}^T\int_{\R^d} |\nabla_x \ln p_t(x)|^2 p_t(x)dxdt\leq \frac{2}{\underline{a}}\bigg(&\int_{\R^d} \ln p_{t_0}(x) p_{t_0}(x)dx+\frac{1}{2}\E[|X_T|^2]+\frac{d}{2}\ln(2\pi)\\&+\int_{t_0}^T\int_{\R^d}\sum_{i=1}^d\partial_{x_i}p_t(x)\bigg(b_i(t,x)-\frac{1}{2}\sum_{j=1}^d\partial_{x_j}a_{ij}(t,x)\bigg)dxdt\bigg).
\end{align*}
When $a\in C^{0,2}_b({\mathcal M}_d(\R))$ and $b\in C^{0,1}_b(\R^d)$ with spatial derivatives of respective orders $2$ and $1$ globally Hölder continuous in space, the Gaussian bounds for $p_t$ and $\nabla_x p_t$ deduced from  Theorems 4.5 and 4.7 in \cite{friedman}, ensure that the estimation \eqref{integderlogrho} should hold for $\rho=2$ as soon as the time integral is restricted to the interval $[t_0,T]$ with $t_0>0$. To our knowledge, even with such a restriction of the time-interval, \eqref{integderlogrho} is not available in the literature for $\rho>2$.\par In fact, we are going to replace the diffusion by another Euler scheme~$\tilde{X}$ with time step $T/M$ and estimate the Wasserstein distance between the marginal laws of the two Euler schemes. We take advantage of the local Gaussian properties of the Euler scheme on each time-step to check that \eqref{integderlogrho} holds when $p_t$ is replaced by $\bar{p}_t$ and to get rid of the boundary terms when performing spatial integration by parts. Finally, we obtain an estimation of the Wasserstein distance between the marginal laws of the diffusion and the Euler scheme by letting $M\to\infty$. Note that we need less spatial regularity on the coefficients $\sigma$ and $b$ than in Theorem 2.2 in \cite{AJK} which directly estimates ${\cal W}_\rho(p_t,\bar{p}_t)$ in dimension $d=1$ by using the optimal coupling given by the inverse transform sampling. 
\begin{aprop}\label{mainprop}
   Under the assumptions of Theorem~\ref{mainthm}, for any $\rho\geq 1$, there exists a finite constant $C$ such that
\begin{equation}\label{majo_intermed}\forall N,M\geq 1,\;\sup_{t\in[0,T]}{\cal W}_\rho({\cal L}(\tilde X_t),{\cal     L}(\bar{X}_t))\leq C \left( \frac{ \left(1+\mathbf{1}_{\gamma=1} \sqrt{\ln(N)}\right)}{N^\gamma}+\frac{ \left(1+\mathbf{1}_{\gamma=1} \sqrt{\ln(M)}\right)}{M^\gamma} \right).
\end{equation}
\end{aprop}
In what follows, we denote the probability density of~$\tilde{X}_t$ for $t\in (0,T]$ by $\tilde{p}_t$ and also set $W_\rho(\tilde{p}_t,\bar{p}_t)=W_\rho({\cal L}(\tilde X_t),{\cal     L}(\bar{X}_t)))$ even for $t=0$ when there is no density.
Let us now explain how we can deduce Theorem \ref{mainthm} from Proposition \ref{mainprop}. Thanks to the triangle inequality, we have
$$ \sup_{t\in[0,T]} W_{\rho}(p_t,\bar{p}_t) \le \sup_{t\in[0,T]} W_{\rho}(p_t,\tilde{p}_t)+\sup_{t\in[0,T]} W_{\rho}(\tilde{p}_t,\bar{p}_t).$$
From the strong error estimate given by Kanagawa~\cite{Ka} in the Lipschitz case and Proposition~14 of Faure~\cite{Faure} for coefficients Hölder continuous in time (see also Theorem 4.1 in Yan~\cite{Yan}), we obtain  $\sup_{t\in[0,T]} W_{\rho}(p_t,\tilde{p}_t)\le \sup_{t\in[0,T]} \E^{1/\rho}[|\tilde{X}_t-X_t|^\rho] \underset{M\rightarrow + \infty}{\rightarrow} 0$, and then deduce Theorem~\ref{mainthm} from~\eqref{majo_intermed}. Note that since the Wasserstein distance is lower semicontinuous with respect to the narrow convergence, the convergence in law of $\tilde{X}_t$ towards~$X_t$ would be enough to obtain the same conclusion.

Concerning the fourth simplifying hypothesis introduced at the beginning of Subsection \ref{Subsec_rig}, we see that the equation~\eqref{der_ags1} given by the results of Ambrosio Gigli and Savaré already gives ``for free'' the first of the two spatial integrations by parts needed to deduce~\eqref{evolwass} from~\eqref{der_Wass}. We will not be able to prove the second integration by parts on the diffusion terms as in~\eqref{evolwass}, but the regularity of the optimal transport maps is sufficient to get an inequality instead of the equality in~\eqref{evolwass} and to go on with the calculations.

The proof is structured as follows. First, we state the optimal transport results between the two Euler schemes $\bar{X}$ and $\tilde{X}$.  Then, we show the Fokker-Planck equation for the Euler scheme and deduce an explicit expression for $\frac{d}{dt}W_{\rho}(\tilde{p}_t,\bar{p}_t)$. Next, we show how we can perform the integration by parts. Last, we put the pieces together and conclude the proof.

\subsection{The optimal transport for the Wasserstein distance $W_{\rho}(\tilde{p}_t,\bar{p}_t)$}\label{subs_rig_opt_tr}
From \eqref{schemeul} and since $\sigma$ does not vanish, it is clear that, for $t>0$, $\bar{X}_t$ and $\tilde{X}_t$ admit positive densities $\bar{p}_t$ and $\tilde{p}_t$ with respect to the Lebesgue measure.

By Theorem 6.2.4 of Ambrosio, Gigli and Savaré~\cite{ags}, for $t\in(0,T]$, there exist measurable optimal transport maps : $\tilde{T}_t,\bar{T}_t:\R^d\to\R^d$ such that $\tilde{T}_t(\tilde{X}_t)$ and $\bar{T}_t(\bar{X}_t)$ have respective densities $\bar{p}_t$ and $\tilde{p}_t$ and \begin{equation}W_\rho^\rho(\tilde{p}_t,\bar{p}_t)=\int_{\R^d}|x-\tilde{T}_t(x)|^\rho \tilde{p}_t(x)dx=\int_{\R^d}|x-\bar{T}_t(x)|^\rho\bar{p}_t(x)dx.\label{egalw}\end{equation}
  Moreover, the positivity of the densities $\tilde{p}_t$ and $\bar{p}_t$, combined with Theorem 3.3.11 and Remark 3.3.14 (b) of Rachev and R\"uschendorf~\cite{raru} ensure that $$dx\;a.e.,\;\tilde{T}_t(x)\in \partial_\rho\tilde{\psi}_t(x)\mbox{ and }\bar{T}_t(x)\in\partial_\rho\bar{\psi}_t(x),$$
where $\tilde{\psi}_t$ and $\bar{\psi}_t$ : $\R^d\to[-\infty,+\infty]$ are two $\rho$-convex (see \eqref{rhoconv}) functions satisfying the duality equation 
\begin{equation}
 \tilde{\psi}_t(x)=-\inf_{y\in\R^d}\left(|x-y|^\rho+\bar{\psi}_t(y)\right) \mbox{ and } \bar{\psi}_t(y)=-\inf_{x\in\R^d}\left(|x-y|^\rho+\tilde{\psi}_t(x)\right).\label{dual2}
\end{equation}
We recall that
\begin{align}
   \partial_\rho \tilde{\psi}_t(x)&=\{y\in\R^d:\tilde{\psi}_t(x)=-(|x-y|^\rho+\bar{\psi}_t(y))\},\label{diff1b}\\
\partial_\rho\bar{\psi}_t(x)&=\{y\in\R^d:\bar{\psi}_t(x)=-(|x-y|^\rho+\tilde{\psi}_t(y))\}.\label{diff2b}
\end{align}
Let us stress that $\bar{T}_t(x)$ now denotes the optimal transport from the law of $\bar{X}_t$ to the law of~$\tilde{X}_t$, while, in Section \ref{sec:characopttransp}, it denoted the optimal transport from the law of $\bar{X}_t$ to the one of $X_t$. However, there is no possible confusion since we will only work in the remainder of Section~\ref{Sec_rig} with the coupling between~$\bar{X}_t$ and~$\tilde{X}_t$.
By the uniqueness in law of the optimal coupling, see e.g Theorem 6.2.4 of Ambrosio, Gigli and Savaré~\cite{ags}, 
  $(\tilde{X}_t,\tilde{T}_t(\tilde{X}_t))$, $(\bar{T}_t(\bar{X}_t),\bar{X}_t)$, $(\bar{T}_t(\bar{X}_t),\tilde{T}_t(\bar{T}_t(\bar{X}_t)))$ and $(\bar{T}_t(\tilde{T}_t(\tilde{X}_t)),\tilde{T}_t(\tilde{X}_t))$ have the same distribution. The equality of the laws of $(\tilde{X}_t,\tilde{T}_t(\tilde{X}_t))$ and $(\bar{T}_t(\tilde{T}_t(\tilde{X}_t)),\tilde{T}_t(\tilde{X}_t))$ implies that $\bar{p}_t(y)dy$ a.e. ${\cal L}(\tilde{X}_t|\tilde{T}_t(\tilde{X}_t)=y)$ and ${\cal L}(\bar{T}_t(\tilde{T}_t(\tilde{X}_t))|\tilde{T}_t(\tilde{X}_t)=y)$ are both equal to the Dirac mass at $\bar{T}_t(y)$ so that $\tilde{X}_t=\bar{T_t}(\tilde{T}_t(\tilde{X}_t))$ a.s.. By positivity of the densities and symmetry we deduce that 
\begin{equation}
   \label{invt}
dx\;a.e.,\;x=\bar{T}_t(\tilde{T}_t(x))=\tilde{T}_t(\bar{T}_t(x)).
\end{equation}

Since, for $\rho \ge 2$, the function $c(x,y)=|x-y|^\rho$ satisfies the conditions (Super), (Twist), (locLip), (locSC) and (H$\infty$) in \cite{villani}, Theorems 10.26-10.28 of Villani~\cite{villani} ensure that $\tilde{\psi}_t$ and $\bar{\psi_t}$ are locally Lipschitz continuous, locally semi-convex, differentiable outside a set of dimension $d-1$, and satisfy
\begin{align}
dx\;a.e.,\;\nabla{\tilde{\psi}}_t(x)+\rho|\tilde{T}_t(x)-x|^{\rho-2}(x-\tilde{T}_t(x))=\nabla\bar{\psi}_t(x)+\rho|\bar{T}_t(x)-x|^{\rho-2}(x-\bar{T}_t(x))=0.\label{gradfaibl}\end{align}
Let us be more precise on the semi-convexity property. When $\rho=2$, we have $\bar{\psi}_t(x)+|x|^2=\sup_{y\in \R^d} \{2 x.y - (\tilde{\psi}_t(y)+|y|^2)\}$ and $\tilde{\psi}_t(x)+|x|^2=\sup_{y\in \R^d}\{2 x.y - (\bar{\psi}_t(y)+|y|^2)\}$, and these functions are convex as they are the suprema of convex functions. When $\rho>2$, we show in Lemma~\ref{lem_semicvx} below that there is a finite constant $C_r$ such that $\bar{\psi}_t(x)+C_r(|x|^2+|x|^\rho)$ and $\tilde{\psi}_t(x)+C_r(|x|^2+|x|^\rho)$ are convex on $B(r)$, where $B(r)=\{x \in \R^d, \ |x|\le r \}$ denotes the ball in~$\R^d$ centered in~$0$ with radius~$r>0$. 

From Theorem 14.25 of Villani~\cite{villani} also known as Alexandrov's second differentiability theorem, we deduce that there is a Borel subset ${\cal A}(\bar{\psi}_t) $ of $\R^d$ such that $\R^d \setminus {\cal A}(\bar{\psi}_t)$ has zero Lebesgue measure and for any $x\in {\cal A}(\bar{\psi}_t)$,  $\bar{\psi}_t$ is differentiable at~$x$ and there is a symmetric matrix~$\nabla^2_A\bar{\psi}_t(x)\in{\cal M}_d(\R)$ called the Hessian of~$\bar{\psi}_t$ such that 
\begin{equation}\label{DL_O2_Alexandrov}\bar{\psi}_t(x+v) \underset{v \rightarrow 0}{=} \bar{\psi}_t(x)+\nabla \bar{\psi}_t(x).v+ \frac{1}{2} \nabla^2_A\bar{\psi}_t(x)v.v +o(|v|^2).  
\end{equation}
Besides, according to Dudley \cite{dudley} p.167, $\nabla^2_A\bar{\psi}_t(x)dx$ coincides with the absolutely continuous part of the distributional Hessian of~$\bar{\psi}_t$, and, by \cite{dudley},  the singular part is positive semidefinite in the following sense : for any $\mathcal{C}^\infty$ function~$\phi$ with compact support on~$\R^d$ with values in the subset of ${\cal M}_d(\R)$ consisting in symmetric positive semidefinite matrices,
\begin{equation}\label{ineg_distr}
\int_{\R^d}\sum_{i,j=1}^d\partial_{x_i}\bar{\psi}_t(x)\partial_{x_j}\phi_{ij}(x)dx \le - \int_{\R^d} \Tr(\nabla^2_A\bar{\psi}_t(x)\phi(x))dx.
\end{equation}

 From~\eqref{DL_O2_Alexandrov}, we can write the second order optimality condition for the minimization of $y\mapsto |x-y|^\rho+\bar{\psi}_t(y)$ and get that
\begin{align*}
   \forall x\in\R^d,\;\forall y\in\partial_\rho\tilde \psi_t(x)\cap{\cal A}(\bar{\psi}_t), \;
\nabla^2_A\bar{\psi}_t(y)+\rho|y-x|^{\rho-2}\left(I_d+(\rho-2)\frac{x-y}{|x-y|}\frac{(x-y)^*}{|x-y|}\right)   \ge 0,
\end{align*}
i.e. it is a positive semidefinite matrix. By Lemma \ref{lem:zeroleb}, \begin{equation}\label{difalexpsibar}
   dx\mbox{ a.e. },\tilde{T}_t(x)\in\partial_\rho\tilde{\psi}_t(x)\cap{\cal A}(\bar{\psi}_t).
\end{equation}

We deduce that
\begin{align}\label{hessfaibl1}
 dx\; a.e.,\;&\nabla^2_A\bar{\psi}_t( \tilde{T}_t(x))+\rho|\tilde{T}_t(x)-x|^{\rho-2}\left(I_d+(\rho-2)\frac{x-\tilde{T}_t(x)}{|x-\tilde{T}_t(x)|}\frac{(x-\tilde{T}_t(x))^*}{|x-\tilde{T}_t(x)|}\right)\geq 0,
\end{align}
and similarly,
\begin{align}
dx\; a.e.,\;&\nabla^2_A{\tilde{\psi}}_t(\bar{T}_t(x))+\rho|\bar{T}_t(x)-x|^{\rho-2}\left(I_d+(\rho-2)\frac{x-\bar{T}_t(x)}{|x-\bar{T}_t(x)|}\frac{(x-\bar{T}_t(x))^*}{|x-\bar{T}_t(x)|}\right)\geq 0.\label{hessfaibl2}
\end{align}

\begin{arem}
One may wonder whether the optimal transport maps $\tilde{T}_t(x)$ and $\bar{T}_t(x)$ satisfy additional regularity properties allowing to proceed as in the heuristic proof, for example to obtain the optimality conditions~\eqref{eul1} and~\eqref{eul2}. We were not able to prove rigorously those conditions. In particular, the assumptions (C) and (STwist) made in Chapter 12 \cite{villani} to get smoothness results are not satisfied by our cost function $c(x,y)=|x-y|^\rho$ for $\rho>2$. Fortunately, the regularity and optimality properties of the optimal transport maps that we have stated from the beginning of Section \ref{subs_rig_opt_tr} will be enough to complete the proof of Theorem \ref{mainthm}. 
\end{arem}
We set \begin{equation}
   \tilde{\tau}_t=\lfloor \frac{Mt}{T}\rfloor\frac{T}{M},\;\tilde{a}(t,x)=\E(a(\tilde{\tau}_t,\tilde{X}_{\tilde\tau_t})|\tilde X_t=x)\mbox{ and }\tilde{b}(t,x)=\E(b(\tilde{\tau}_t,\tilde{X}_{\tilde\tau_t})|\tilde X_t=x).\label{defttaut}
\end{equation} The rest of Section~\ref{Sec_rig} will consist in proving the following result.
\begin{aprop}\label{prop_majo_wasser} Let us suppose that $$\exists K\in[0,+\infty), \forall x \in \R^d, \ \sup_{t\in[0,T]} |\sigma(t,x)|+|b(t,x)| \le K(1+|x|)$$ and assume uniform ellipticity : there exists a positive constant $\underline{a}$ such that $a(t,x)-\underline{a} I_d$ is positive semidefinite for any $(t,x)\in[0,T]\times\R^d$. Then, $t\mapsto W^\rho_\rho(\tilde{p}_t,\bar{p}_t)$ is absolutely continuous and such that $dt$ a.e.,
\begin{align*}
   \frac{d}{dt}W_\rho^\rho(\tilde{p}_t,\bar{p}_t)\le &   C \Bigg( W_\rho^\rho(\tilde{p}_t,\bar{p}_t) +\int_{\R^d}|x-\bar{T}_t(x)|^{\rho-1} |\bar{b}(t,x)-b(t,x)| \bar{p}_t(x)dx
 \\ &+\int_{\R^d}|x-\tilde{T}_t(x)|^{\rho-1} |\tilde{b}(t,x)-b(t,x)| \tilde{p}_t(x)dx\\
&+\int_{\R^d} |x-\bar{T}_t(x)|^{\rho-2} \Tr[   (  \bar{a}(t,x)-a(t,x) )^2] \bar{p}_{t}(x)dx  \\
&+\int_{\R^d} |x-\tilde{T}_t(x)|^{\rho-2} \Tr[   (  a(t,x) -\tilde{a}(t,x) )^2] \tilde{p}_{t}(x)dx\Bigg),\end{align*}
where the finite constant $C$ does not depend on $t\in[0,T]$, $x_0 \in \R^d$ and $N,M\geq 1$.
\end{aprop}
With this result, we can repeat the arguments of Subsection~\ref{subs_Gronwall}, and obtain Proposition~\ref{mainprop} and thus Theorem~\ref{mainthm}.


\subsection{Proof of Proposition~\ref{prop_majo_wasser} }

The proof is based on the second of the two next propositions which estimates the time-derivative of the Wasserstein distance under gradually stronger assumptions on the coefficients~$a$ and $b$.  

\begin{aprop}\label{prop_1}
We assume ellipticity : $a(t,x)$ is positive definite for any $t\in(0,T]$, $x\in \R^d$. We also suppose that $\exists K\in[0,+\infty), \forall x \in \R^d, \ \sup_{t\in[0,T]} |\sigma(t,x)|+|b(t,x)| \le K(1+|x|)$. Then $t\mapsto W^\rho_\rho(\tilde{p}_t,\bar{p}_t)$ is absolutely continuous and such that $dt$ a.e.,
   \begin{align}
   \frac{d}{dt}W_\rho^\rho(\tilde{p}_t,\bar{p}_t)\le &   -\frac{1}{2} \int_{\R^d}\Tr[\nabla^2_A\tilde{\psi}_t(x)\tilde{a}(t,x)] \tilde{p}_{t}(x)dx \notag \\& -\frac{1}{2} \int_{\R^d}\Tr[\nabla^2_A\bar{\psi}_t(\tilde{T}_t(x))\bar{a}(t,\tilde{T}_t(x))] \tilde{p}_{t}(x)dx \label{der_ags3}\\&+\rho\int_{\R^d}|\tilde{T}_t(x)-x|^{\rho-2}(\tilde{T}_t(x)-x).\left(\bar{b}(t,\tilde{T}_t(x))-\tilde{b}(t,x)\right)\tilde{p}_t(x)dx.\notag\end{align}
\end{aprop}

\begin{aprop}\label{prop_2}
Under the assumptions of Proposition \ref{prop_majo_wasser}, $dt$ a.e.,
\begin{align}
   \frac{d}{dt}W_\rho^\rho(\tilde{p}_t,\bar{p}_t)\le &   \frac{ \rho(\rho-1)^2 }{8 \underline{a} }  \int_{\R^d}  |\tilde{T}_t(x)-x|^{\rho-2} \Tr[   (  \bar{a}(t,\tilde{T}_t(x))-\tilde{a}(t,x))^2] \tilde{p}_{t}(x)dx \notag \\&+\rho\int_{\R^d}|\tilde{T}_t(x)-x|^{\rho-2}(\tilde{T}_t(x)-x).\left(\bar{b}(t,\tilde{T}_t(x))-\tilde{b}(t,x)\right)  \tilde{p}_t(x)dx. \label{result_majo2}\end{align}
\end{aprop}
\begin{arem}
   Notice that these two propositions still hold with $$\tilde{a}(t,x)=\E(\hat{\sigma}\hat{\sigma}^*(\tilde{\tau}_t,\tilde{X}_{\tilde\tau_t})|\tilde X_t=x)\mbox{ and }\tilde{b}(t,x)=\E(\hat{b}(\tilde{\tau}_t,\tilde{X}_{\tilde\tau_t})|\tilde X_t=x)$$ when $\tilde{X}_t$ is the Euler scheme with step $T/M$ for the stochastic differential equation 
$$Y_t=y_0+\int_0^t\hat{b}(s,Y_s)ds+\int_0^t\hat{\sigma}(s,Y_s)dW_s,\;t\leq T$$
with $y_0\in\R^d$, $\hat{b}:[0,T]\times\R^d\to \R^d$ and $\hat{\sigma} : [0,T]\times\R^d\to {\cal M}_d(\R)$ satisfying the same conditions as $b$ and $\sigma$.
\end{arem}

Proposition \ref{prop_majo_wasser} is deduced from Proposition \ref{prop_2} by using the triangle inequalities
\begin{align*}
|\bar{b}(t,\tilde{T}_t(x))-\tilde{b}(t,x)| &\le  
|\bar{b}(t,\tilde{T}_t(x))-b(t,\tilde{T}_t(x))|+|b(t,\tilde{T}_t(x))-b(t,x)| +|b(t,x)-\tilde{b}(t,x)|,\\
\frac{1}{3} \Tr[   (  \bar{a}(t,\tilde{T}_t(x))-\tilde{a}(t,x))^2] & \le \Tr[   (  \bar{a}(t,\tilde{T}_t(x)) - a(t,\tilde{T}_t(x)) )^2]+ \Tr[   ( a(t,\tilde{T}_t(x))-a(t,x))^2]\\
&+\Tr[   ( a(t,x)- \tilde{a}(t,x))^2],
\end{align*}
the bounds on the first derivatives of~$a$ and $b$ and $\tilde{T}_t\#\tilde p_t=\bar p_t$.

 The proofs of Propositions~\ref{prop_1} and \ref{prop_2} are given in the two next sections.

\subsubsection{Proof of Proposition \ref{prop_1}}\label{subsec_prop44}
The proof of Proposition~\ref{prop_1} is split in the next three paragraphs. We first explicit the time evolution of the probability density of the Euler scheme. Then, this enables us to apply the results of Ambrosio, Gigli and Savaré and get a formula for  $\frac{d}{dt} W_\rho^\rho(\tilde{p}_t,\bar{p}_t)$ in \eqref{der_ags2}. Last, we show that we have the desired inequality by a spatial integration by parts. Of course, we work under the assumptions of Proposition~\ref{prop_1} in these two paragraphs.

\paragraph{The Fokker-Planck equation for the Euler scheme.}

We focus on the Euler scheme~$\bar{X}$ and use the notations given in the introduction.

For $k\in\{0,\hdots,N\}$, denoting by $\bar{\mu}_{t_k}$ the law of $\bar{X}_{t_k}$, one has that for $t\in(t_k,t_{k+1}]$, the law of $(\bar{X}_{t_k},\bar{X}_t)$ is $\bar{\mu}_{t_k}(dy)G^{a,b}_{t_k,t}(y,x)dx$ where
\begin{align*}
   G^{a,b}_{t_k,t}(y,x)=\frac{e^{-\frac{1}{2(t-t_k)}(x-y-b(t_k,y)(t-t_k)).a^{-1}(t_k,y)(x-y-b(t_k,y)(t-t_k))}}{(2\pi(t-t_k))^{d/2}\sqrt{\det(a(t_k,y))}}.
\end{align*}
Notice that $\bar{\mu}_0(dy)=\delta_{x_0}(dy)$ while for $k\geq 1$, $\bar{\mu}_{t_k}(dy)=\bar{p}_{t_k}(y)dy$.
\begin{alem}\label{lem:regufpt} The function
   \begin{equation}
   \bar{v}_t(x)=\bar{b}(t,x)-\frac{1}{2\bar{p}_t(x)}\int_{\R^d}a(\tau_t,y)\nabla_xG^{a,b}_{\tau_t,t}(y,x)\bar{\mu}_{\tau_t}(dy)\label{decvbar}
\end{equation}
defined for $t\in[0,T]\setminus\{t_0,t_1,\hdots,t_N\}$ and $x\in\R^d$
is such that $\int_0^T \left( \int_{\R^d} |\bar{v}_t(x)|^\rho \bar{p}_t(x)dx\right)^{1/\rho}dt<\infty$ and $\partial_t \bar{p}_t + \nabla .(\bar{v}_t \bar{p}_t)=0$ holds in the sense of distributions on $(0,T)\times\R^d$.
\end{alem}

\begin{adem}
   Let $\varphi$ be a $C^\infty$ function with compact support on $(0,T) \times \R^d$. From~\eqref{eul_conti}, we apply Ito's formula to $\varphi(t,\bar{X}_t)$ between $0$ and~$T$ and then take the expectation to get
\begin{align*}
0&=  \int_0^T \E\left[ \partial_t \varphi(t,\bar{X}_t)+\nabla_x \varphi(t,\bar{X}_t).b(\tau_t,\bar{X}_{\tau_t})+\frac{1}{2}\Tr\left(\nabla^2_x \varphi(t,\bar{X}_t) a(\tau_t,\bar{X}_{\tau_t})\right) \right]dt \\
&=  \int_0^T \E\left[ \partial_t \varphi(t,\bar{X}_t)+ \nabla_x \varphi(t,\bar{X}_t).\E[ b(\tau_t,\bar{X}_{\tau_t})|\bar{X}_t]+\frac{1}{2}\Tr\left(\nabla_x^2 \varphi(t,\bar{X}_t) a(\tau_t,\bar{X}_{\tau_t})\right) \right]dt ,
\end{align*}
from the tower property of the conditional expectation. This then leads to:
\begin{align*}
  0&=\int_{0}^{T}\left[\int_{\R^d}(\partial_t \varphi(t,x)+\bar{b}(t,x).\nabla_x\varphi(t,x))\bar{p}_t(x)+\frac{1}{2}\int_{\R^d}\Tr(a(\tau_t,y)\nabla^2_x\varphi(t,x))G^{a,b}_{\tau_t,t}(y,x)\bar{\mu}_{\tau_t}(dy)\right]dxdt.\end{align*}
By performing one integration by parts with respect to~$x$, we get that
\begin{equation}
   \partial_t\bar{p}_t(x)+\nabla.(\bar{v}_t(x)\bar{p}_t(x))=0\label{transpobar}
\end{equation} holds in the sense of distributions in $(0,T)\times\R^d$. \\It remains to check that 
$\int_0^T \left( \int_{\R^d} |\bar{v}_t(x)|^\rho \bar{p}_t(x)dx\right)^{1/\rho}dt<\infty$. From the assumption on~$b$ and $\sigma$, the Euler scheme has bounded moments, and therefore
\begin{align*}
 \int_0^T \left( \int_{\R^d} |\bar{b}(t,x)|^\rho \bar{p}_t(x)dx\right) &=\int_0^T \E[|\E(b(t,\bar{X}_{\tau_t})|\bar{X}_t) |^\rho] dt \\
&\le \int_0^T K^\rho 2^{\rho-1}(1+\E[|\bar{X}_{\tau_t}|^\rho])dt < \infty.
\end{align*}

We can then focus on the second term in \eqref{decvbar}. 
 We notice that for $t\in (t_k,t_{k+1})$, we have
\begin{align*}
   \left|\frac{1}{\bar{p}_t(x)}\int_{\R^d}a(t_k,y)\nabla_xG^{a,b}_{t_k,t}(y,x)\bar{\mu}_{t_k}(dy)\right|^\rho&=\frac{1}{(t-t_k)^\rho}\left|\frac{\int_{\R^d}(x-y-b(t_k,y)(t-t_k))G^{a,b}_{t_k,t}(y,x)\bar{\mu}_{t_k}(dy)}{\int_{\R^d}G^{a,b}_{t-t_k}(y,x)\bar{\mu}_{t_k}(dy)}\right|^\rho\\
&\leq \frac{\int_{\R^d}|x-y-b(t_k,y)(t-t_k)|^\rho G^{a,b}_{t_k,t}(y,x)\bar{\mu}_{t_k}(dy)}{(t-t_k)^\rho\bar{p}_{t}(x)},
\end{align*}
by Jensen's inequality and using $\bar{p}_{t}(x)=\int_{\R^d}G^{a,b}_{t-t_k}(y,x)\bar{\mu}_{t_k}(dy)$.\\Since $$G^{a,b}_{t_k,t}(y,x)=2^{d/2} G^{2a,b}_{t_k,t}(y,x)e^{-\frac{1}{4(t-t_k)}(x-y-b(t_k,y)(t-t_k)).a^{-1}(t_k,y)(x-y-b(t_k,y)(t-t_k))}$$ and $\max_{z\ge 0}z^{\rho/2}e^{-\alpha z}=\left(\frac{\rho}{2\alpha e} \right)^{\rho/2}$ for $\alpha>0$, we get
$$|x-y-b(t_k,y)(t-t_k)|^\rho G^{a,b}_{t_k,t}(y,x)\leq 2^{d/2} \left(\frac{2\rho \bar{\lambda}(a(t_k,y)) (t-t_k)}{e}\right)^{\rho/2}G^{2a,b}_{t_k,t}(y,x),$$ 
where $\bar{\lambda}(a)$ denotes the largest eigenvalue of the matrix $a$. Therefore, 
 $$ \left|\frac{1}{\bar{p}_t(x)}\int_{\R^d}a(t_k,y)\nabla_xG^{a,b}_{t_k,t}(y,x)\bar{\mu}_{t_k}(dy)\right|^\rho  \leq \frac{ 2^{d/2}}{\bar{p}_t(x)} \int_{\R^d} \left(\frac{2\rho K(1+|y|)^2 }{e(t-t_k)}\right)^{\rho/2}G^{2a,b}_{t_k,t}(y,x) \bar{\mu}_{t_k}(dy),$$
since by assumption $\bar{\lambda}(a(t,x))\le K(1+|x|)^2$ for some $K<+\infty$,  and
we deduce that 
\begin{equation}
   \left( \int_{\R^d}\left|\frac{1}{\bar{p}_t(x)}\int_{\R^d}a(t_k,y)\nabla_xG^{a,b}_{t_k,t}(y,x)\bar{\mu}_{t_k}(dy)\right|^\rho\bar{p}_t(x)dx \right)^{1/\rho}\leq 2^{\frac{d}{2 \rho}} \left(\frac{2\rho K}{e(t-t_k)}\right)^{1/2}   \E[(1+|\bar{X}_{t_k}|)^{\rho}]^{1/\rho}.\label{contvlrho}
\end{equation}
Using  $\int_0^T(t-\tau_t)^{-1/2}dt=2\sqrt{NT}$ and the boundedness of the moments of the Euler scheme, we get that $\int_0^T \left( \int_{\R^d} |\bar{v}_t(x)|^\rho \bar{p}_t(x)dx\right)^{1/\rho}dt<\infty$. 
\end{adem}

\paragraph{The time derivative of the Wasserstein distance.}

To compute $\frac{d}{dt} W_\rho^\rho(\tilde{p}_t,\bar{p}_t)$, we are going to adapt to the differentiation of the Wasserstein distance between two absolutely continuous curves the proof of Theorem~8.4.7 of Ambrosio, Gigli and Savaré~\cite{ags} where one of these curves is constant. We also need to introduce \begin{equation*}
   \tilde{v}_t(x)=\tilde{b}(t,x)-\frac{1}{2\tilde{p}_t(x)}\int_{\R^d}a(\tilde\tau_t,y)\nabla_xG^{a,b}_{\tilde\tau_t,t}(y,x)\tilde{\mu}_{\tilde\tau_t}(dy)
\end{equation*}
where $\tilde\tau_t$ is defined in \eqref{defttaut} and $\tilde{\mu}_{\tilde\tau_t}(dy)$ denotes the law of $\tilde{X}_{\tilde\tau_t}$. 
Note that the conclusion of Lemma \ref{lem:regufpt} is also valid with $(\bar{p}_t,\bar{v}_t)$ replaced by $(\tilde{p}_t,\tilde{v}_t)$.
By the last statement in Theorem~8.3.1 \cite{ags},
$t\mapsto \bar{p}_t(x)dx$ and $t\mapsto \tilde{p}_t(x)dx$ are absolutely continuous curves in the set of probability measures on $\R^d$ with bounded moment of order $\rho$ endowed with $W_\rho$ as a metric. By the triangle inequality, one deduces that $t\mapsto W_\rho(\bar{p}_t,\tilde{p}_t)$ is an absolutely continuous function, which, with the continuous differentiability of $w\mapsto|w|^\rho$ on $\R$, ensures the absolute continuity of $t\mapsto W^\rho_\rho(\bar{p}_t,\tilde{p}_t)$. By the first statement in Theorem~8.3.1 and Proposition~8.4.6 \cite{ags}, there exist Borel vector fields $\bar{v}_t^\star(x)$ and $\tilde{v}^\star_t(x)$ defined on $[0,T]\times \R^d$ satisfying $\int_0^T \left( \int_{\R^d} |\bar{v}^\star_t(x)|^\rho \bar{p}_t(x)dx\right)^{1/\rho}+\left( \int_{\R^d} |\tilde{v}^\star_t(x)|^\rho \tilde{p}_t(x)dx\right)^{1/\rho}dt<\infty$, \begin{equation}
   \partial_t \bar{p}_t + \nabla .(\bar{v}^\star_t \bar{p}_t)=0=\partial_t \tilde{p}_t + \nabla .(\tilde{v}^\star_t \tilde{p}_t)\mbox{ in the sense of distributions on $(0,T)\times\R^d$},\label{eqtt}
\end{equation} and $dt$ a.e. on $(0,T)$, 
$$\lim_{h\to 0}\frac{W_\rho(\bar{p}_{t+h},(\mathbf{i}+h\bar{v}^\star_t)\#\bar{p}_t)+W_\rho(\tilde{p}_{t+h},(\mathbf{i}+h\tilde{v}^\star_t)\#\tilde{p}_t)}{h}=0,$$
where $\mathbf{i}(x)=x$ denotes the identity function on $\R^d$. Note that these vector fields characterized (up to $dt$ a.e. equality) by \eqref{eqtt} together with
$$dt\mbox{ a.e. }\left( \int_{\R^d} |\bar{v}^\star_t(x)|^\rho \bar{p}_t(x)dx\right)^{1/\rho}\leq \lim_{s\to t}\frac{W_\rho(\bar{p}_s,\bar{p}_t)}{|s-t|}\mbox{ and }\left( \int_{\R^d} |\tilde{v}^\star_t(x)|^\rho \tilde{p}_t(x)dx\right)^{1/\rho}\leq \lim_{s\to t}\frac{W_\rho(\tilde{p}_s,\tilde{p}_t)}{|s-t|}$$
are called in Proposition 8.4.5 \cite{ags} the tangent vectors to the absolutely continuous curves $t\mapsto \bar{p}_t(x)dx$ and $t\mapsto \tilde{p}_t(x)dx$.

Since, by the triangle inequality,
\begin{align*}
&W_\rho(\bar{p}_{t+h},\tilde{p}_{t+h})\le W_\rho(\bar{p}_{t+h},(\mathbf{i}+h\bar{v}^\star_t)\#\bar{p}_t)+W_\rho((\mathbf{i}+h\bar{v}^\star_t)\#\bar{p}_t,(\mathbf{i}+h\tilde{v}^\star_t)\#\tilde{p}_t) +W_\rho((\mathbf{i}+h\tilde{v}^\star_t)\#\tilde{p}_t,\tilde{p}_{t+h}) \\
&W_\rho((\mathbf{i}+h\bar{v}^\star_t)\#\bar{p}_t,(\mathbf{i}+h\tilde{v}^\star_t)\#\tilde{p}_t) \le W_\rho((\mathbf{i}+h\bar{v}^\star_t)\#\bar{p}_t,\bar{p}_{t+h})+W_\rho(\bar{p}_{t+h},\tilde{p}_{t+h})+W_\rho(\tilde{p}_{t+h},(\mathbf{i}+h\tilde{v}^\star_t)\#\tilde{p}_t),
\end{align*}
one deduces that $dt$ a.e.,
$$\frac{d}{dt}W_\rho^\rho(\bar{p}_t,\tilde{p}_t)=\lim_{h\rightarrow 0} \frac{W_\rho^\rho((\mathbf{i}+hv_t)\#\bar{p}_t,(\mathbf{i}+h\tilde{v}^\star_t)\#\tilde{p}_t)-W_\rho^\rho(\bar{p}_t,\tilde{p}_t)}{h}.$$
Since $(\mathbf{i}+h\bar{v}^\star_t)(\bar{X}_t)$ and $(\mathbf{i}+h\tilde{v}^\star_t)(\bar{T}_t(X_t))$ are respectively distributed according to $(\mathbf{i}+h\bar{v}^\star_t)\#\bar{p}_t$ and $(\mathbf{i}+h\tilde{v}^\star_t)\#\tilde{p}_t$ and $\int_{\R^d}|x-\bar{T}_t(x)|^\rho \bar{p}_t(x)dx=W_\rho^\rho(\bar{p}_t,\tilde{p}_t)$,
\begin{align*}
&W_\rho^\rho((\mathbf{i}+h\bar{v}^\star_t)\#\bar{p}_t,(\mathbf{i}+h\tilde{v}^\star_t)\#\tilde{p}_t)\le \int_{\R^d}|x+h\bar{v}^\star_t(x)-\bar{T}_t(x)-h \tilde{v}^\star_t(\bar{T}_t(x))|^\rho \bar{p}_t(x)dx \\
&= W_\rho^\rho(\bar{p}_t,\tilde{p}_t)+h \int_{\R^d} \rho |x-\bar{T}_t(x)|^{\rho-2} (x-\bar{T}_t(x)).(\bar{v}^\star_t(x)-\tilde{v}^\star_t(\bar{T}_t(x)))\bar{p}_t(x)dx +o(h).
\end{align*}
Letting $h\to 0^+$ and $h\to 0^-$, one respectively deduces that the two next inequalities hold $dt$ a.e.
\begin{align*}
   \rho \int_{\R^d} |x-\bar{T}_t(x)|^{\rho-2} (x-\bar{T}_t(x))&.(\bar{v}^\star_t(x)-\tilde{v}^\star_t(\bar{T}_t(x)))\bar{p}_t(x)dx\geq \frac{d}{dt}W_\rho^\rho(\bar{p}_t,\tilde{p}_t)\\
&\geq \rho\int_{\R^d}  |x-\bar{T}_t(x)|^{\rho-2} (x-\bar{T}_t(x)).(\bar{v}^\star_t(x)-\tilde{v}^\star_t(\bar{T}_t(x)))\bar{p}_t(x)dx.
\end{align*}

Like in Remark 8.4.8 \cite{ags}, Proposition 8.5.4 in the same book \cite{ags} justifies that $(\bar{v}_t^\star,\tilde{v}_t^\star)$ may be replaced by $(\bar{v}_t,\tilde{v}_t)$ in this formula. With \eqref{invt} and $\tilde{p}_t=\bar{T}_t\#\bar{p}_t$, we deduce that, $dt$ a.e, 
\begin{equation*}
\frac{d}{dt} W_\rho^\rho(\tilde{p}_t,\bar{p}_t)=\rho\int_{\R^d}|\tilde{T}_t(x)-x|^{\rho-2}(x-\tilde{T}_t(x)).\tilde{v}_t(x)\tilde{p}_t(x)+|\bar{T}_t(x)-x|^{\rho-2}(x-\bar{T}_t(x)).\bar{v}_t(x)\bar{p}_t(x)dx.
\end{equation*}

Using~\eqref{gradfaibl}, plugging the expressions of $\bar{v}_t$ and $\tilde{v}_t$ then \eqref{invt} and $\bar{T}_t\#\bar{p}_t=\tilde{p}_t$, we get that, $dt$ a.e.,
\begin{align}
   \frac{d}{dt}W_\rho^\rho(\tilde{p}_t,\bar{p}_t)=&-\int_{\R^d}\nabla \tilde{\psi}_t(x).\tilde{v}_t(x)\tilde{p}_t(x)+\nabla\bar{\psi}_t(x).\bar{v}_t(x)\bar{p}_t(x)dx\notag\\=&\frac{1}{2}\int_{\R^d}\int_{\R^d}\nabla\tilde{\psi}_t(x).a(\tilde{\tau}_t,y)\nabla_x G^{a,b}_{\tilde{\tau}_t,t}(y,x)dx\tilde{\mu}_{\tau_t}(dy)\notag \\&+ \frac{1}{2}\int_{\R^d}\int_{\R^d}\nabla\bar{\psi}_t(x).a(\tau_t,y)\nabla_x G^{a,b}_{\tau_t,t}(y,x)dx\bar{\mu}_{\tau_t}(dy) \label{der_ags2}\\&+\rho\int_{\R^d}|\tilde{T}_t(x)-x|^{\rho-2}(\tilde{T}_t(x)-x).\left(\bar{b}(t,\tilde{T}_t(x))-\tilde{b}(t,x)\right)\tilde{p}_t(x)dx.\notag\end{align}

\paragraph{The integration by parts inequality.}

The aim of this paragraph is to prove the following inequality
\begin{equation}\label{Ineg_IPP}
\int_{\R^d}\int_{\R^d}\nabla\bar{\psi}_t(x).a(\tau_t,y)\nabla_x G^{a,b}_{\tau_t,t}(y,x)dx\bar{\mu}_{\tau_t}(dy)  \le -\int_{\R^d}\Tr[\nabla^2_A\bar{\psi}_t(x)\bar{a}(t,x)] \bar{p}_{t}(x)dx.
\end{equation}
To do so, we introduce cutoff functions to use the inequality~\eqref{ineg_distr}. We recall that $B(r)$ denotes the closed ball in~$\R^d$ centered in~$0$ with radius~$r>0$. For  ${\ell}\geq 1$, we consider a $C^\infty$ function $\varphi_{\ell}:\R^d\to [0,1]$ such that:
\begin{equation*}
\forall x \in B({\ell}), \varphi_{\ell}(x)=1, \  \forall x \not \in B(2{\ell}), \varphi_{\ell}(x)=0 \text{ and } \forall x \in \R^d, |\nabla\varphi_{\ell}(x)|\leq \frac{2}{{\ell}}.
\end{equation*}
One has
\begin{align*}
   \int_{\R^d}\nabla\bar{\psi}_t(x).a(\tau_t,y)&\nabla_x G^{a,b}_{\tau_t,t}(y,x)dx=\int_{\R^d}\nabla\bar{\psi}_t(x).a(\tau_t,y)\nabla_x (\varphi_{\ell}(x)G^{a,b}_{\tau_t,t}(y,x))dx\\
&+\int_{\R^d}\nabla\bar{\psi}_t(x).a(\tau_t,y)\left((1-\varphi_{\ell}(x))\nabla_x G^{a,b}_{\tau_t,t}(y,x)-G^{a,b}_{\tau_t,t}(y,x)\nabla\varphi_{\ell}(x)\right)dx.
\end{align*}

From~\eqref{gradfaibl} and \eqref{egalw}, we have $\int_{\R^d}|\nabla\bar{\psi}_t(x)|^{\frac{\rho}{\rho-1}}\bar{p}_t(x)dx=\rho^{\frac{\rho}{\rho-1}}W_\rho^\rho(\tilde{p}_t,\bar{p}_t)$. By~\eqref{contvlrho} and H\"older's inequality, we deduce that 
\begin{align*}
\int_{\R^d}|\nabla\bar{\psi}_t(x)|&\left|\int_{\R^d}a(\tau_t,y)\nabla_xG^{a,b}_{\tau_t,t}(y,x)\bar{\mu}_{\tau_t}(dy)\right|dx \\
&\leq   2^{\frac{d}{2 \rho}} \left(\frac{2\rho K}{e(t-t_k)}\right)^{1/2}   \E[(1+|\bar{X}_{t_k}|)^{\rho}]^{1/\rho} \times \rho W_\rho^{\rho-1}(\tilde{p}_t,\bar{p}_t).
\end{align*}
We also have 
\begin{align*}
\int_{\R^d\times\R^d}|\nabla\bar{\psi}_t(x).a(\tau_t,y)\nabla\varphi_{\ell}(x)| G^{a,b}_{\tau_t,t}(y,x)\bar{\mu}_{\tau_t}(dy)dx&\leq  \frac{2}{\ell} \int_{\R^d\times\R^d}|\nabla\bar{\psi}_t(x)| \bar{\lambda}(a(\tau_t,y)) G^{a,b}_{\tau_t,t}(y,x)\bar{\mu}_{\tau_t}(dy)dx.\\
& \leq    \frac{K}{\ell}  \E[(1+|\bar{X}_{\tau_t}|)^{2\rho}]^{1/\rho}\times  \rho W_\rho^{\rho-1}(\tilde{p}_t,\bar{p}_t). 
\end{align*}
Using the dominated convergence theorem, we obtain
$$\lim_{{\ell}\to\infty}\int_{\R^d\times\R^d}\nabla\bar{\psi}_t(x).a(\tau_t,y)\left((1-\varphi_{\ell}(x))\nabla_x G^{a,b}_{\tau_t,t}(y,x)-G^{a,b}_{\tau_t,t}(y,x)\nabla\varphi_{\ell}(x)\right)\bar{\mu}_{\tau_t}(dy)dx=0.$$
On the other hand we use the inequality~\eqref{ineg_distr} to get 
\begin{align*}
   \int_{\R^d} \nabla\bar{\psi}_t(x).a(\tau_t,y)\nabla_x(\varphi_{\ell}(x)G^{a,b}_{\tau_t,t}(y,x)) dx \leq -\int_{\R^d}\Tr(\nabla^2_A\bar{\psi}_t(x)a(\tau_t,y))\varphi_{\ell}(x)G^{a,b}_{\tau_t,t}(y,x)dx,
\end{align*}
for any $y\in \R^d$, and thus
\begin{align}
  \int_{\R^d\times\R^d}\nabla\bar{\psi}_t(x)&.a(\tau_t,y)\nabla_x G^{a,b}_{\tau_t,t}(y,x)dx\bar{\mu}_{\tau_t}(dy)\notag\\
&\leq -\limsup_{{\ell}\to\infty} \int_{\R^d\times\R^d}\Tr(\nabla^2_A\bar{\psi}_t(x)a(\tau_t,y))\varphi_{\ell}(x)G^{a,b}_{\tau_t,t}(y,x)\bar{\mu}_{\tau_t}(dy)dx\notag\\
&=-\limsup_{{\ell}\to\infty} \int_{\R^d}\Tr(\nabla^2_A\bar{\psi}_t(x)\bar{a}(t,x))\varphi_{\ell}(x)\bar{p}_t(x)dx,\label{limsupell}
\end{align}
where we used the definition of $\bar{a}$ for the equality. Using this definition again, we get\begin{align*}
\int_{\R^d} |x-\bar{T}_t(x)|^{\rho-2}|\bar{a}(t,x)|\bar{p}_t(x)dx &= \int_{\R^d\times\R^d} |x-\bar{T}_t(x)|^{\rho-2}|a(\tau_t,y)| G^{a,b}_{\tau_t,t}(y,x)dx\bar{\mu}_{\tau_t}(dy) \\
&\le W_\rho^{\rho-2}(\tilde p_t,\bar{p}_t) \left( \int_{\R^d} (K(1+|y|)^2)^{\rho/2}\bar{\mu}_{\tau_t}(dy) \right)^{2/\rho}<\infty.
\end{align*}
With \eqref{hessfaibl1}, we deduce that $\Tr(\nabla^2_A\bar{\psi}_t(x)\bar{a}(t,x))\bar{p}_{t}(x)$ is the sum of a non-negative and an integrable function. Using Fatou's Lemma for the contribution of the non-negative function and Lebesgue's theorem for the contribution of the integrable function in \eqref{limsupell}, we finally obtain~\eqref{Ineg_IPP}.\\By symmetry, we have
$$\int_{\R^d}\int_{\R^d}\nabla\tilde{\psi}_t(x).a(\tilde{\tau}_t,y)\nabla_x G^{a,b}_{\tilde{\tau}_t,t}(y,x)dx\tilde{\mu}_{\tau_t}(dy)\leq -\int_{\R^d}\Tr[\nabla^2_A\tilde{\psi}_t(x)\tilde{a}(t,x)] \tilde{p}_{t}(x)dx.$$
Using $\tilde{T}_t\#\tilde{p}_t=\bar{p}_t$ in the right-hand-side of \eqref{Ineg_IPP} leads to
$$\int_{\R^d}\int_{\R^d}\nabla\bar{\psi}_t(x).a(\tau_t,y)\nabla_x G^{a,b}_{\tau_t,t}(y,x)dx\bar{\mu}_{\tau_t}(dy)  \le -\int_{\R^d}\Tr[\nabla^2_A\bar{\psi}_t(\tilde{T}_t(x))\bar{a}(t,\tilde{T}_t(x))] \tilde{p}_{t}(x)dx.$$
Plugging the two last inequalities in~\eqref{der_ags2} gives Proposition~\ref{prop_1}.

\subsubsection{Proof of Proposition~\ref{prop_2}}

Let $h(x)=|x|^\rho$. We have $\nabla h(x)=\rho|x|^{\rho-2}x$ and $(\nabla h)^{-1}(x)=1_{\{x\neq 0\}}\rho^{-\frac{1}{\rho-1}}|x|^{\frac{2-\rho}{\rho-1}}x$. Notice that when $\rho=2$, $(\nabla h)^{-1}(x)=\frac{x}{2}$ is also defined for $x=0$.  

By~\eqref{gradfaibl}, we have $dx$ a.e. $\bar{T}_t(x)=x+(\nabla h)^{-1}(\nabla\bar{\psi}_t(x))$, $\tilde{T}_t(x)=x+(\nabla h)^{-1}(\nabla{\tilde{\psi}}_t(x))$. Using~\eqref{invt} and Lemma \ref{lem:zeroleb} with ${\cal A}=\{x\in\R^d:\bar{T}_t(x)=x+(\nabla h)^{-1}(\nabla\bar{\psi}_t(x))\}$, we deduce that
$$dx\;a.e.,\;x=x+(\nabla h)^{-1}(\nabla{\tilde{\psi}}_t(x))+(\nabla h)^{-1}(\nabla\bar\psi_t(x+(\nabla h)^{-1}(\nabla{\tilde{\psi}}_t(x)))),$$
and thus
\begin{equation}\label{lien_gradients}
dx\ a.e.,\; \nabla{\tilde{\psi}}_t(x)=-\nabla\bar\psi_t(x+(\nabla h)^{-1}(\nabla{\tilde{\psi}}_t(x))).
\end{equation}
When $\rho=2$, $\nabla^*(\nabla h)^{-1}(x)=\frac{1}{2}I_d$ and when $\rho>2$, $\nabla^*(\nabla h)^{-1}(x)=\rho^{-\frac{1}{\rho-1}} |x|^{\frac{2-\rho}{\rho-1}} \left( I_d+\frac{2-\rho}{\rho-1} \frac{xx^*}{|x|^2}\right)$ for $x\not =0$. Because of the singularity of $\nabla^*(\nabla h)^{-1}(x)$ at the origin for $\rho>2$, we set $\mathcal{E}=\{x\in \R^d, \tilde{T}_t(x) \not =x \}$ if $\rho>2$ and $\mathcal{E}=\R^d$ if $\rho=2$.\\
By \eqref{difalexpsibar}, Lemma \ref{lem_semicvx} and Property (i) in Theorem 14.25 of Villani~\cite{villani}, we can thus perform first order expansions in equation~\eqref{lien_gradients} to get that $dx \text{ a.e. on } \mathcal{E}$,
\begin{equation}
   \nabla^2_A{\tilde{\psi}}_t(x)=-\nabla^2_A\bar{\psi}_t(x+(\nabla h)^{-1}(\nabla{\tilde{\psi}}_t(x)))\left[I_d+\nabla^*(\nabla h)^{-1}(\nabla{\tilde{\psi}}_t(x))\nabla^2_A{\tilde{\psi}}_t(x)\right].\label{lienhessal}
\end{equation}
Using~\eqref{gradfaibl}, we get
$$\nabla^*(\nabla h)^{-1}(\nabla{\tilde{\psi}}_t(x) )=\frac{1}{\rho} |x-\tilde{T}_t(x)|^{2-\rho}\left( I_d+\frac{2-\rho}{\rho-1}  v_x v_x^*\right), \;dx \text{ a.e. on } \mathcal{E}$$
with $v_x=\frac{x-\tilde{T}_t(x)}{|x-\tilde{T}_t(x)|}$. We define the positive definite matrix $A(x)=I_d+(\rho-2)v_x v_x^*$ with inverse $A^{-1}(x)=I_d+\frac{2-\rho}{\rho-1}  v_x v_x^*$. Plugging the above identities in \eqref{lienhessal}, we obtain
\begin{equation}\label{lien_hessiennes}
\nabla^2_A{\tilde{\psi}}_t(x)+\nabla^2_A\bar{\psi}_t(\tilde{T}_t(x))=-\frac{1}{\rho}|x-\tilde{T}_t(x)|^{2-\rho}\nabla^2_A\bar{\psi}_t(\tilde{T}_t(x))A^{-1}(x)\nabla^2_A{\tilde{\psi}}_t(x), \ dx \text{ a.e. on } \mathcal{E} .
\end{equation}
 We set $M(x)=\frac{1}{\rho}|x-\tilde{T}_t(x)|^{2-\rho}\nabla^2_A{\tilde{\psi}}_t(x) + A(x)$ for $x\in \mathcal{E}$ such that the right-hand-side makes sense. By \eqref{hessfaibl1}, Lemma \ref{lem:zeroleb} and \eqref{invt}, $M(x)$ is a positive semidefinite matrix $dx \text{ a.e. on } \mathcal{E}$. Moreover,
$$\nabla^2_A{\tilde{\psi}}_t(x)=\rho|x-\tilde{T}_t(x)|^{\rho-2}(M(x)-A(x)), \ dx \text{ a.e. on } \mathcal{E} .$$ Using this equality in the right hand side of~\eqref{lien_hessiennes}, we get
$ \nabla^2_A{\tilde{\psi}}_t(x)=-\nabla^2_A\bar{\psi}_t(\tilde{T}_t(x))A^{-1}(x)M(x) $, which gives
$$-\nabla^2_A\bar{\psi}_t(\tilde{T}_t(x))A^{-1}(x)M(x)=\rho|x-\tilde{T}_t(x)|^{\rho-2}(M(x)-A(x)),  \ dx \text{ a.e. on } \mathcal{E}. $$
Therefore $dx$ a.e. on $\mathcal{E}$, every element of $\R^d$ in the kernel of the matrix $M(x)$ belongs to the kernel of the invertible matrix $A(x)$ so that $M(x)$ is invertible. We finally have
$$-\nabla^2_A\bar{\psi}_t(\tilde{T}_t(x))=\rho|x-\tilde{T}_t(x)|^{\rho-2}(A(x)-A(x)M^{-1}(x)A(x)),  \ dx \text{ a.e. on } \mathcal{E}.$$
Plugging this equality in \eqref{der_ags3}, we obtain that 
\begin{align}
&   \frac{d}{dt}W_\rho^\rho(\tilde{p}_t,\bar{p}_t)  \le \rho\int_{\R^d}|\tilde{T}_t(x)-x|^{\rho-2}(\tilde{T}_t(x)-x).\left(\bar{b}(t,\tilde{T}_t(x))-\tilde{b}(t,x)\right)\tilde{p}_t(x)dx\notag\\
 & +  \frac{1}{2} \int_{\mathcal{E}} \rho |x-\tilde{T}_t(x)|^{\rho-2} \Tr[   (A(x)-M(x)) \tilde{a}(t,x)+(A(x)-A(x)M^{-1}(x)A(x)) \bar{a}(t,\tilde{T}_t(x))] \tilde{p}_{t}(x)dx\notag\\
&-\frac{1}{2} \int_{\R^d\setminus\mathcal{E} }\Tr[\nabla^2_A\tilde{\psi}_t(x)\tilde{a}(t,x)+\nabla^2_A\bar{\psi}_t(\tilde{T}_t(x))\bar{a}(t,\tilde{T}_t(x))] \tilde{p}_{t}(x)dx .\label{der_ags6} \end{align}

When $\rho>2$ and $x \not \in \mathcal{E}$, we have from~\eqref{invt}, \eqref{hessfaibl1},~\eqref{hessfaibl2} and Lemma \ref{lem:zeroleb} that $\nabla^2_A\tilde{\psi}_t(x)$ and $\nabla^2_A\bar{\psi}_t(\tilde{T}_t(x))$ are positive semidefinite $dx$ a.e. on $\R^d\setminus \mathcal{E}$ and therefore
$$\Tr[\nabla^2_A\tilde{\psi}_t(x)\tilde{a}(t,x)+\nabla^2_A\bar{\psi}_t(\tilde{T}_t(x))\bar{a}(t,\tilde{T}_t(x))] \ge 0  \ dx \text{ a.e. on }\R^d\setminus \mathcal{E}. $$
Therefore the third term in the right-hand-side of \eqref{der_ags6} is non positive. Using Lemma~\ref{lemma_Wp2} for the second term, we conclude that~\eqref{result_majo2} holds by remarking that the definition of ${\mathcal E}$ ensures that
$$\int_{\R^d\setminus{\mathcal E}}  |\tilde{T}_t(x)-x|^{\rho-2} \Tr[   (  \bar{a}(t,\tilde{T}_t(x))-\tilde{a}(t,x))^2] \tilde{p}_{t}(x)dx=0.$$

\section{Technical Lemmas}\label{sectchlem}

\subsection{Transport of negligible sets}
\begin{alem}\label{lem:zeroleb}
   Let $\bar{T}(x)$ and $\tilde{T}(x)$ be measurable optimal transport maps for $W_\rho$ with $\rho\geq 2$ between two probability measures with positive densities $\bar{p}$ and $\tilde{p}$ with respect to the Lebesgue measure on $\R^d$ : $\bar{p}=\tilde{T}\#\tilde{p}$ and $\tilde{p}=\bar{T}\#\bar{p}$. For any Borel subset ${\cal A}$ of $\R^d$ such that $\R^d\setminus{\cal A}$ has zero Lebesgue measure, $dx$ a.e. $\tilde{T}(x)\in{\cal A}$ and $\bar{T}(x)\in{\cal A}$.
\end{alem}
\begin{dem}
Since $\tilde{T}\#\tilde{p}=\bar{p}$ and $\R^d\setminus{\cal A}$ has zero Lebesgue measure, $$\int_{\R^d}1_{\cal A}(\tilde{T}(x))\tilde{p}(x)dx=\int_{\R^d}1_{\cal A}(x)\bar{p}(x)dx=1.$$By positivity of $\tilde{p}$, one concludes that $dx$ a.e. $\tilde{T}(x)\in{\cal A}$.
\end{dem}

\subsection{A key Lemma on pseudo-distances between matrices}
The next Lemma holds as soon as $\rho>1$ and not only under the assumption $\rho\geq 2$ made from Section \ref{sec:characopttransp} on.
\begin{alem}\label{lemma_Wp2}
For $v\in \R^d$ such that $|v|=1$, let $A$ denote the positive definite matrix $I_d +(\rho-2)v v^*$. Let $M,a_1,a_2\in{\cal M}_d(\R)$ be positive definite symmetric matrices.  Then for any $\underline{a}>0$ such that $a_i-\underline{a}I_d$ is positive semidefinite for $i\in\{1,2 \}$, one has \begin{equation}
   \Tr \left[A\left\{(I_d-A^{-1}M)a_1+(I_d-M^{-1}A)a_2 \right\} \right]\leq  \frac{(1\vee (\rho-1))^2 }{4 \underline{a}(1 \wedge (\rho-1)) }  \Tr\left[\left(a_1-a_2 \right)^2 \right].\label{maot}
\end{equation}
\end{alem}
Notice that the left-hand side of the inequality is linear in $a_1$ and $a_2$, whereas thanks to the positivity of $\underline{a}$ we obtain the quadratic factor $\Tr\left[\left(a_1-a_2 \right)^2 \right]$ in the right-hand side.
\begin{dem}
We define $\tilde{M}=A^{- \frac{1}{2}} M A^{- \frac{1}{2}}$, where $A^{- \frac{1}{2}}$ is the inverse of the square-root $A^{\frac{1}{2}}$ of the symmetric positive definite matrix~$A$. Let $\mathbf{T}=\Tr \left[A\left\{(I_d-A^{-1}M)a_1+(I_d-M^{-1}A)a_2 \right\} \right]$ denote the quantity to be estimated. We have, using the cyclicity of the trace for the third equality below, 
\begin{align*}
\mathbf{T} &=\Tr \left[(A-M)a_1+(A-AM^{-1}A)a_2  \right] \\
&=\Tr \left[A^{\frac{1}{2}} (I_d-\tilde{M})A^{\frac{1}{2}} a_1+A^{\frac{1}{2}}(I_d-\tilde{M}^{-1})A^{\frac{1}{2}}a_2  \right] \\
&= \Tr \left[ (I_d-\tilde{M})\left\{ A^{\frac{1}{2}} a_1 A^{\frac{1}{2}}- \tilde{M}^{-1} A^{\frac{1}{2}}a_2 A^{\frac{1}{2}}  \right\} \right].
\end{align*}
Let $(\lambda_1,\hdots,\lambda_d)$ denote the vector of eigenvalues of the symmetric positive definite matrix $\tilde{M}$, $D(\lambda_1,\hdots,\lambda_d)$ be the diagonal matrix with diagonal coefficients $\lambda_1,\hdots,\lambda_d$ and $O$ be the orthogonal matrix such that $\tilde{M}=O^*D(\lambda_1,\hdots,\lambda_d)O$. We define \begin{align*}
   (I_d\vee \tilde{M})^{-1}&:=O^*D((1\vee \lambda_1)^{-1},\hdots,(1\vee \lambda_d)^{-1})O,\\
(I_d-\tilde{M})^+&:=O^*D((1-\lambda_1)^{+},\hdots,(1-\lambda_d)^{+})O,\\(\tilde{M}-I_d)^+&:=O^*D((\lambda_1-1)^{+},\hdots,(\lambda_d-1)^{+})O.
\end{align*} 
Since for all $\lambda\in\R$, $1-\lambda=(1-\lambda)(1\vee \lambda)^{-1}-\lambda^{-1}((\lambda-1)^+)^2$ and $(1-\lambda)\lambda^{-1}=(1-\lambda)(1\vee \lambda)^{-1}+\lambda^{-1}((1-\lambda)^+)^2$, we have
\begin{align}
\mathbf{T}=&\Tr\bigg[(I_d-\tilde{M})(I_d\vee \tilde{M})^{-1} [ A^{\frac{1}{2}}(a_1-a_2)A^{\frac{1}{2}} ]\bigg]\notag\\&-\Tr\bigg[\tilde{M}^{-1}((\tilde{M}-I_d)^+)^2 A^{\frac{1}{2}}a_1A^{\frac{1}{2}}\bigg]-\Tr\bigg[\tilde{M}^{-1}((I_d-\tilde{M})^+)^2 A^{\frac{1}{2}}a_2 A^{\frac{1}{2}}\bigg].\label{majot}
\end{align}
On the one hand, by Cauchy-Schwarz and Young's inequalities, for symmetric matrices $S_1,S_2$,
$$\Tr(S_1S_2)\le \sqrt{\Tr(S_1^2)\Tr(S_2^2)}\leq \underline{a}(1 \wedge (\rho-1))\Tr(S_1^2)+ \frac{1}{4 \underline{a}(1 \wedge (\rho-1)) }\Tr(S_2^2),$$ 
which implies that
\begin{align*}&\Tr\bigg[(I_d-\tilde{M})(I_d\vee \tilde{M})^{-1} [ A^{\frac{1}{2}}(a_1-a_2)A^{\frac{1}{2}} ]\bigg] \\ &\le \underline{a}(1 \wedge (\rho-1)) \sum_{i=1}^d \frac{(1-\lambda_i)^2}{(1\vee \lambda_i)^2}+ \frac{1}{4 \underline{a}(1 \wedge (\rho-1)) } \Tr\bigg[\left(A^{\frac{1}{2}}(a_1-a_2)A^{\frac{1}{2}}\right)^2  \bigg].
\end{align*}
On the other hand, we recall that $\Tr(S_1S_2)\ge c \Tr(S_1)$ when $S_1,S_2$ are symmetric positive semidefinite matrices such that $S_2-cI_d$ is  positive semidefinite. Since the smallest eigenvalue of $A$ is $1\wedge(\rho-1)$,  $A^{\frac{1}{2}}a_1A^{\frac{1}{2}}-\underline{a}(1 \wedge (\rho-1))I_d$ is positive semidefinite and we get
$$\Tr\bigg[\tilde{M}^{-1}((\tilde{M}-I_d)^+)^2 A^{\frac{1}{2}}a_1A^{\frac{1}{2}}\bigg]\ge  \underline{a}(1 \wedge (\rho-1))\sum_{i=1}^d\frac{((\lambda_i-1)^+)^2}{\lambda_i}, $$ and similarly 
 $$\Tr\bigg[\tilde{M}^{-1}((I_d-\tilde{M})^+)^2 A^{\frac{1}{2}}a_2 A^{\frac{1}{2}}\bigg]\ge  \underline{a}(1 \wedge (\rho-1))\sum_{i=1}^d\frac{((1-\lambda_i)^+)^2}{\lambda_i}. $$
Since $\frac{(1-\lambda_i)^2}{(1\vee \lambda_i)^2}-\frac{((\lambda_i-1)^+)^2}{\lambda_i}-\frac{((1-\lambda_i)^+)^2}{\lambda_i} \le 0$, we finally get that:
\begin{align*}
 \mathbf{T} &\le  \frac{1}{4 \underline{a}(1 \wedge (\rho-1)) }\Tr\bigg[\left(A^{\frac{1}{2}}(a_1-a_2)A^{\frac{1}{2}}\right)^2 \bigg]\\
&\le \frac{(1\vee (\rho-1))^2 }{4 \underline{a}(1 \wedge (\rho-1)) }\Tr\left[\left(a_1-a_2\right)^2 \right].
\end{align*}
We have used for the last inequality the cyclicity of the trace and $\Tr(AS)\le(1 \vee (\rho-1))\Tr(S) $ for any positive semidefinite matrix~$S$, since the largest eigenvalue of~$A$ is $1 \vee (\rho-1)$.
\end{dem}
\begin{arem}
\begin{enumerate}
   \item In dimension~$d=1$, the only eigenvalue of~$A$ is $\rho-1$, and we get the slightly better bound $$  A\left\{(1-A^{-1}M)a_1+(1-M^{-1}A)a_2 \right\} \leq  \frac{ (\rho-1) }{4 \underline{a}}  \left(a_1-a_2 \right)^2  .$$
\item Inequality \eqref{maot} still holds with $\Tr((a_1-a_2)^2)$ replaced by $\Tr((a_1-a_2)(a_1-a_2)^*]$ in the right-hand side for all $a_1,a_2\in{\cal M}_d(\R)$ such that $a_1+a_1^*-2\underline{a}I_d$ and $a_2+a_2^*-2\underline{a}I_d$ are positive semidefinite.
\item Since the second and third terms in the right-hand-side of \eqref{majot} are non-positive, applying Cauchy-Schwarz inequality to the first term, one obtains that $\forall a_1,a_2\in{\cal M}_d(\R)$,
 $$\Tr \left[A\left\{(I_d-A^{-1}M)a_1+(I_d-M^{-1}A)a_2 \right\} \right]\leq \sqrt{(d+\rho-2)(1\vee(\rho-1))}\sqrt{\Tr((a_1-a_2)(a_1-a_2)^*)}.$$

\end{enumerate}
   
\end{arem}

\subsection{Semi-convexity of $\rho$-convex functions for $\rho>2$}

\begin{alem}\label{lem_semicvx}
 Let $\rho>2$ and $t\in (0,T]$. Under the framework of Subsection~\ref{subs_rig_opt_tr}, for any $r\in(0,+\infty)$, there is a finite constant $C_r$ such that $x\mapsto \bar{\psi}_t(x)+C_r(|x|^2+|x|^\rho)$ and $x\mapsto\tilde{\psi}_t(x)+C_r(|x|^2+|x|^\rho)$ are convex on the closed ball $B(r)$ centered at the origin with radius $r$.
\end{alem}
\begin{dem}
We do the proof for $\tilde{\psi}_t$ and follow the arguments of Figalli and Gigli~\cite{FigalliGigli}. Let $r\in(0,+\infty)$. We consider the set
$$A=\{y \in \R^d, \exists x \in B(r), \ \tilde{\psi}_t(x)\le -|x-y|^\rho-\bar{\psi}_t(y)+1 \}.$$ 
Let us check that the existence of a finite constant $K_{r,\rho}$ depending on $r$ and $\rho$ such that $\sup_{y\in A} \min_{x\in B(r)}|x-y| \le K_{r,\rho}$ ensures that the conclusion holds. We have $A\subset B(K_{r,\rho}')$ with $K'_{r,\rho}=K_{r,\rho}+r$.  This gives that
$$\forall x \in B(r), \ \tilde{\psi}_t(x) =\sup_{y\in A} -(\bar{\psi}_t(y)+|x-y|^\rho)=\sup_{y \in B(K'_{r,\rho}) }   -(\bar{\psi}_t(y)+|x-y|^\rho). $$
We also remark that for a constant $C_r$ large enough, $x\mapsto -|x-y|^\rho +C_r(|x|^2+|x|^\rho)$ is convex for any $y\in B(K'_{r,\rho})$. In fact, the Hessian matrix
$$ - \rho|x-y|^{\rho-2}\left(I_d + (\rho-2) \frac{(x-y)(x-y)^*}{|x-y|^2}\right) + C_r\left( 2I_d + \rho |x|^{\rho-2}\left(I_d + (\rho-2) \frac{xx^*}{|x|^2}\right)\right) $$
is positive semidefinite for $C_r$ large enough since for any $y\in B(K'_{r,\rho})$ and $x \in \R^d$, $|x-y|^{\rho-2}\le 2^{(\rho-3)^+}(({K'_{r,\rho}})^{\rho-2}+|x|^{\rho-2})$. Thus, for $x\in B(r)$,  $\tilde{\psi}_t(x)+C_r(|x|^2+|x|^\rho)$ is convex as it is the supremum of convex functions. 

We now prove that $\sup_{y\in A} \min_{x\in B(r)}|x-y| \le K_{r,\rho}$.
Let $y\in A$. If $y\in B(r+1)$, we have $\min_{x\in B(r)}|x-y|\le 1$.  When $|y|>r+1$, we consider $x \in B(r)$ such that $\tilde{\psi}_t(x)\le -|x-y|^\rho-\bar{\psi}_t(y)+1$. We have for $x'\in \R^d$,
\begin{align*}
\tilde{\psi}_t(x') \ge -\bar{\psi}_t(y)-|x'-y|^\rho &= -\bar{\psi}_t(y)-|x-y|^\rho+|x-y|^\rho-|x'-y|^\rho \\
&\ge \tilde{\psi}_t(x)-1+|x-y|^\rho-|x'-y|^\rho.
\end{align*}
We have $|x-y|\ge 1$ and we  take $x'=x - \lambda(x-y)$ with $\lambda\in [0,1/|x-y|]$ so that $|x'|\le r+1$. We get
$$ \tilde{\psi}_t(x')- \tilde{\psi}_t(x)+1 \ge  |x-y|^\rho (1-(1-\lambda)^\rho) .$$
There is $\eta\in (0,1)$ such that 
$ \forall \lambda \in [0,\eta],  1-(1-\lambda)^\rho \ge \frac{\rho}{2} \lambda $. We choose  $\lambda=\eta /|x-y|$  and get
$$  \tilde{\psi}_t(x')- \tilde{\psi}_t(x)+1 \ge \frac{\rho}{2} \eta |x-y|^{\rho-1}.$$
Therefore $|x-y|\le \left( \frac{2}{\rho \eta}  [\sup_{x'\in B(r+1)}\tilde{\psi}_t(x')- \inf_{x\in B(r)}\tilde{\psi}_t(x)+1 ]\right)^{1/(\rho-1)}$ with the function $\tilde{\psi}_t$ locally bounded since it is locally Lipschitz according to Theorem 10.26 \cite{villani}.
\end{dem}

\subsection{Estimations using Malliavin calculus}
\begin{alem}\label{malcal}Under the assumptions of Theorem \ref{mainthm},
  we have for all $\rho\geq 1$ :
\[
\exists C<+\infty,\;\forall N\geq 1,\;\forall t\in[0,T],\;\mathbb{E}\left[  \left\vert \mathbb{E}\left[  W_{t}-W_{\tau_{t}%
}|\bar{X}_{t}\right]  \right\vert ^{\rho}\right]  \leq
C\left((t-\tau_t)\wedge \left (\frac{(t-\tau_t)^2}t+\frac{1}{N^2}\right )\right)^{\rho/2}.
\]
\end{alem}
\begin{adem}[of Lemma \ref{malcal}]
By Jensen's inequality, 
\begin{align*}
  \E\left[|\E(W_t-W_{\tau_t}|\bar{X}_t)|^\rho\right]&\leq \E\left[|W_t-W_{\tau_t}|^{\rho}\right]\leq C(t-\tau_t)^{\rho/2}.
\end{align*}
Let us now check that the left-hand-side is also smaller than $C\left (\frac{(t-\tau_t)^2}t+\frac{1}{N^2}\right )^{\rho/2}$.
To do this, we will study
\[
{\mathbb{E}}\left[  \langle W_{t}-W_{\tau_{t}},g(\bar{X}_{t})\rangle \right],
\]
where $g:\mathbb{R}^d \rightarrow \mathbb{R}^d$ is any smooth function. 

In order to continue, we need to do various estimations on the
Euler scheme, its limit and their Malliavin derivatives, which we denote by $D^i_u\bar{X}^j_t$ and $D^i_uX^j_t$. Let $\eta
_{t}=\min\{t_{i};t\leq t_{i}\}$ denote the discretization time just after
$t$. We have $D^i_u\bar{X}^j_t=0$ for $u>t$, $i,j=1,...,d$ and 
for $u\leq t$,
\begin{align*}
  D^i_{u}&\bar{X}^j_{t}=1_{\{t\leq
  \eta_u\}}\sigma_{ji}(\tau_t,\bar{X}_{\tau_{t}})\\&+
1_{\{t>\eta_u\}}\sum_{k=1}^d\left(1_{\{k=j\}}+\left (\partial_{x_k}\sigma_{jl}(\tau_t,\bar{X}_{\tau_{t}})(W^l_t-W^l_{\tau_t})+\partial_{x_k}b_j(\tau_t,\bar{X}_{\tau_{t}})(t-\tau_t)\right )\right )D^i_u\bar{X}^k_{\tau_{t}}
. 
\end{align*}
Let us define $D\bar{X}:=(D^i\bar{X}^j)_{ij}$. Then by induction, one clearly obtains that for $u\le t$,
\begin{align}
D_{u}\bar{X}_{t} &  =\sigma(\tau_u,\bar{X}_{\tau_{u}})^*\mathcal{\bar{E}}_{u,t},\label{eq:DbX}\\
\sigma &=(\sigma_{ij})_{ij}\nonumber\\
\nonumber
\mathcal{\bar{E}}_{u,t} &  =\left\{
\begin{array}
[c]{cccc}%
I & \text{if} & \tau_{t}< \eta_{u}\\
\left(  I+{\nabla b}(\tau_t,\bar{X}_{\tau_{t}})(t-\tau_{t})+\sigma^{\prime}(\tau_t,\bar
{X}_{\tau_{t}})(W_{t}-W_{\tau_{t}})\right)   & \text{if} & \eta_{u}=\tau_{t}\\
\prod_{i=\frac{N\eta_{u}}{T}}^{\frac{N\tau_{t}}{T}-1}\left(  I+{\nabla b}(t_i,\bar{X}_{t_{i}%
})(t_{i+1}-t_{i})+\sigma^{\prime}(t_i,\bar{X}_{t_{i}})(W_{t_{i+1}}-W_{t_{i}%
})\right)    & \text{if} &
\eta_{u}<\tau_{t}\\
\times\left(  I+{\nabla b}(\tau_t,\bar{X}_{\tau_{t}})(t-\tau_{t})+\sigma
^{\prime}(\tau_t,\bar{X}_{\tau_{t}})(W_{t}-W_{\tau_{t}})\right) 
.
\end{array}
\right.
\end{align}
Here ${\nabla b}:=(\partial_{x_k}b_j)_{kj}$, $\sigma'=(\partial_{x_k}\sigma_{j\cdot})_{kj}$ and $\prod_{i=1}^nA_i:=A_1\cdots A_n$. Therefore the above product between $\sigma'$ and the increment of $W$ is to be interpreted as the inner product between vectors once $k$ and $j$ are fixed.

Note that $\mathcal{\bar{E}}$ satisfies the following properties:
1. $\mathcal{\bar{E}}_{u,t}=\mathcal{\bar{E}}_{\eta(u),t}$ and
2. $\mathcal{\bar{E}}_{t_i,t_j}\mathcal{\bar{E}}_{t_j,t}=\mathcal{\bar{E}}_{t_i,t}$
for $t_i\le t_j\le t$. 

We also introduce
the  process $\mathcal{E}$ as the $d\times d$-matrix solution to the linear stochastic differential equation
\begin{equation}
\label{eq:E}
\mathcal{E}_{u,t}=I+\int_u^t\mathcal{E}_{u,s}{\nabla b}(s,X_s)ds+\int_u^t \mathcal{E}_{u,s}\sigma'(s,X_s)dW_s.
\end{equation}

The next lemma, the proof of which is postponed at the end of the present
proof states some useful properties of the processes $\mathcal{E}$ and $\bar{\mathcal{E}}$. From now on, for $A\in{\cal M}_d(\R)$, $|A|=\sqrt{\mathrm{Tr}(A^*A)}$ denotes its Frobenius norm. 
\begin{alem}\label{lemme_majorations} Let us assume that $b,\sigma \in C^2_b$. Then, we have:
\begin{eqnarray}
\label{eq:propE}
&&
\sup_{0\leq s\leq t \le T}{\mathbb{E}}\left[
| \mathcal{E}_{s,t}^{-1}|^\rho\right]+{\mathbb{E}}\left[
 |\mathcal{E}_{s,t}|^{\rho}\right]  \leq C,  \
\sup_{0\leq s\leq t \le T}{\mathbb{E}}\left[
 |\mathcal{\bar{E}}_{s,t}|^{\rho}\right]  \leq C, \ \ \ \ \ \\
&&
\sup_{0\leq s,u\leq t \le T}{\mathbb{E}}\left[
| D_u\bar{\mathcal{E}}_{s,t}|^\rho+| D_u \mathcal{E}_{s,t}|^\rho\right ]   \leq
C,\label{eq:A13} \\
&& \sup_{0\leq t \le T}{\mathbb{E}}\left[  \left\vert \mathcal{E}_{0,t}%
-\mathcal{\bar{E}}_{0,t}\right\vert ^{\rho}\right]     \leq \frac{C}{N^{\rho(\frac{1}{2}\wedge \gamma)}}, \label{vitesse_forte}
\end{eqnarray}
where $C$ is a positive constant depending only on $\rho$ and $T$.
\end{alem}

We next define the localization given by
\[
\psi=\varphi\left(  |\mathcal{E}_{0,t}^{-1}\left(  \mathcal{E}_{0,t}%
-\mathcal{\bar{E}}_{0,t}\right )|^2  \right)  .
\]

Here $\varphi:\mathbb{R\rightarrow}[0,1]$ is a~$C^\infty$ symmetric function so
that
\[
\varphi(x)=\left\{
\begin{array}
[c]{ccc}%
0, & \text{if} & |x|>\frac{1}{2},\\
1, & \text{if } & |x|<\frac{1}{4}.
\end{array}
\right.
\]
Note that for $M$ in the open ball $B(I_d,2^{-1/2})$ centered at $I_d$ with radius $2^{-1/2}$, one has that $|M-I_d|<2^{-1/2}$ and therefore the sum $\sum_{j=0}^\infty (I_d-M)^k$ converges absolutely. In other words, the map $M\mapsto M^{-1}$ is well defined and bounded on $B(I_d,2^{-1/2})$. 

Now, as $\varphi(x)=0$ for $|x|>2^{-1}$, then if $\psi>0$  we have that $M:=\mathcal{E}_{0,t}^{-1}
\mathcal{\bar{E}}_{0,t}\in B(I_d,2^{-1/2}) $. Therefore $
\mathcal{\bar{E}}_{0,t}^{-1}$ exists and 
\begin{equation}
\label{eq:23}
|\mathcal{\bar{E}}_{0,t}^{-1}|\le |(\mathcal{{E}}_{0,t}^{-1}\mathcal{\bar{E}}_{0,t})^{-1}||\mathcal{{E}}_{0,t}^{-1}|\le \sum_{k=0}^\infty \frac 1{\sqrt{2}^k}|\mathcal{{E}}_{0,t}^{-1}|.
\end{equation}
One has
\begin{align*}
   \mathbb{E}\left[  \langle W_{t}-W_{\tau_{t}},g(\bar{X}_{t})\rangle \right]
&=\mathbb{E}\left[  \langle W_{t}-W_{\tau_{t}},g(\bar{X}_{t})\rangle\psi\right]+\mathbb{E}\left[  \langle W_{t}-W_{\tau_{t}},g(\bar{X}_{t})\rangle(1-\psi)\right]\\
&=\int_{\tau_t}^t\E\left[\psi \mathrm{Tr}(D_u\bar{X}_{t}\nabla g(\bar{X}_{t}))\right] du+\E\left[\int_{\tau_t}^t\langle D_u\psi,
g(\bar{X}_{t})\rangle  du\right]\\
&+\mathbb{E}\left[  \langle W_{t}-W_{\tau_{t}},g(\bar{X}_{t})\rangle(1-\psi)\right].
\end{align*}
The second equality follows from the duality formula (see e.g. Definition 1.3.1 in \cite{N}). Since for $\tau_{t}\leq u\leq t$%
\begin{align*}
&{\mathbb{E}}\left[  \psi \mathrm{Tr}(D_{u}\bar{X}_{t}\nabla g(\bar{X}_{t}))\right]   
={\mathbb{E}}\left[  \psi \mathrm{Tr}(\sigma(\tau_t, \bar{X}_{\tau_{t}%
})^*\nabla g(\bar{X}_{t})
)\right]\\
&=t^{-1}\E\left[\int_0^t \psi\mathrm{Tr}(\sigma(\tau_t, \bar{X}_{\tau_{t}})^*(D_s\bar{X}_t)^{-1}D_sg(\bar{X}_{t})
)
ds\right]  \\
&  =t^{-1}{\mathbb{E}}\left[  g(\bar{X}_{t})^*\int_{0}^{t} \psi  \sigma(\tau_t, \bar
{X}_{\tau_{t}})^*\mathcal{\bar{E}%
}_{s,t}^{-1}\sigma^{-1}\left( \tau_s,  \bar{X}_{\tau_{s}}\right)^*  \delta W_{s} \right].
\end{align*}
Here $\delta W$ denotes the Skorohod vector integral (see \cite{N}). Then 
one deduces
\begin{align}
\mathbb{E}\left[  \left.  W_{t}-W_{\tau_{t}}\right|  \bar{X}_{t}\right]
&=t^{-1}\int_{\tau_{t}}^{t}\mathbb{E}\left[  \left.  \int_{0}^{t}\psi
\sigma(\tau_t, \bar{X}_{\tau_{t}})^*
\mathcal{\bar{E}}_{s,t}^{-1}\sigma^{-1}\left( \tau_s,  \bar{X}_{\tau_{s}}\right)^*\delta W_{s}\right|  \bar{X}_{t}\right]
du
\notag\\
&+
\mathbb{E}
\left[\int_{\tau_t}^t D_u\psi du \Bigg|  \bar{X}_{t} \right ]+\mathbb{E}\left[  \left.  \left(  W_{t}-W_{\tau_{t}}\right)  (1-\psi
)\right|  \bar{X}_{t}\right]  .\label{espcond}
\end{align}
 In order to obtain the
conclusion of the Lemma, we need to bound the $L^\rho$-norm of each term on the right-hand-side of (\ref{espcond}). In particular, we will use the following estimate (which also proves the existence of the Skorohod integral on the left side below) which can be found in Proposition 1.5.4 in \cite{N}: 
\begin{equation}\label{controle_normep}
\left\Vert \int_{0}^{t}\psi\sigma(\tau_t, \bar{X}_{\tau_{t}})^* \mathcal{\bar{E}}_{s,t}^{-1}\sigma^{-1}\left(\tau_s,
\bar{X}_{\tau_{s}}\right) ^*\delta W_{s}\right\Vert
_{\rho}\leq C(\rho)\left\Vert \psi\sigma(\tau_t, \bar{X}_{\tau_{t}})^*\mathcal{\bar{E}}_{\cdot,t}^{-1}
\sigma^{-1}\left (\tau_{\cdot}, \bar
{X}_{\tau_{\cdot}}\right )^*\right\Vert _{1,\rho},
\end{equation}
where $\|F_\cdot \|_{1,\rho}^\rho =\E\left[\left(\int_0^t| F_s|^2 ds
  \right)^{\rho/2}+\left(\int_0^t \int_0^t |D_uF_s|^2 dsdu
  \right)^{\rho/2}\right]$. By Jensen's inequality for $\rho\ge 2$, we have
\begin{equation}\label{upper_bound_1p} \|F_\cdot \|_{1,p}^\rho \le  t^{\rho/2-1} \int_0^t \E[|F_s|^\rho] ds + t^{\rho-2}
\int_0^t\int_0^t \E[|D_uF_s|^\rho]dsdu,
\end{equation}
and we will use this inequality to upper bound~\eqref{controle_normep}. When
$1\le \rho\le 2$, we will use alternatively the following upper bound $\|F_\cdot
\|_{1,\rho}^\rho \le \left( \int_0^t \E[|F_s|^2] ds
 \right)^{\rho/2}+\left(\int_0^t \int_0^t \E[|D_uF_s|^2] dsdu
  \right)^{\rho/2}$ that comes
from Jensen's inequality.

Note that for any two invertible matrices $A$, $B$ in ${\cal M}_d(\R)$, we have that 
$|B^* A(B^{-1})^*|\le |B||A||B^{-1}|$. Choosing $B=\sigma(\tau_t, \bar{X}_{\tau_{t}}
)$ and $A=\mathcal{\bar{E}}_{s,t}
^{-1}$, remarking that $|B^{-1}|=\sqrt{\Tr(a^{-1}(\tau_t, \bar{X}_{\tau_{t}}
))}$ and using the boundedness of $a$ and the uniform ellipticity, we deduce that there exists a finite constant $ C$ such that 
\begin{eqnarray}
  \int_{0}^{t}{\mathbb{E}}\left[ \left(  \psi|\sigma(\tau_t, \bar{X}_{\tau_{t}}
)^* \mathcal{\bar{E}}_{s,t}
^{-1}\sigma^{-1}\left( \tau_s,  \bar{X}_{\tau_{s}}\right) ^*|\right)  ^{\rho}\right]ds  &\leq& C\int_{0}^{t}{\mathbb{E}}\left[  \psi
^{\rho}|\mathcal{\bar{E}}_{0,t}^{-1}|^\rho|\mathcal{\bar{E}}_{0,\eta (s)}|^{\rho}\right]ds \label{majosigxisig-1}\\ &\leq&C
\sqrt{\E[|\mathcal{{E}}_{0,t}^{-1}|^{2\rho}]}
\int_0^t \sqrt{\E[|\mathcal{\bar{E}}_{0,\eta (s)}|^{2\rho}]}ds \leq C t,\notag
\end{eqnarray}
by using the estimates~\eqref{eq:propE}. Note that we have used that $\psi\mathcal{\bar{E}}^{-1}_{s,t}=\psi
\mathcal{\bar{E}}^{-1}_{0,t}\mathcal{\bar{E}}_{0,\eta(s)}$ and \eqref{eq:23}.

Next, we focus on getting an upper bound for
\begin{equation}
  \int_0^t \int_{0}^{t}{\mathbb{E}}\left[ \left|  D_{u}\left(  \psi\sigma(\tau_t, \bar{X}%
_{\tau_{t}})^* \mathcal{\bar{E}%
}_{s,t}^{-1}\sigma^{-1}\left( \tau_s,  \bar{X}_{\tau_{s}}\right) ^*\right)  \right|^{\rho}\right] ds du  .\label{eq:Dloc}%
\end{equation}
To do so, we compute the above derivative using basic derivation rules, which gives for $l=1,...,d$
\begin{align}
 D^l_{u}\left(  \psi\sigma(\tau_t, \bar{X}_{\tau_{t}})^* \mathcal{\bar{E}}_{s,t}^{-1}\sigma^{-1}\left( \tau_s,  \bar{X}%
_{\tau_{s}}\right)^*\right) & =D^l_u\psi
\sigma(\tau_t, \bar{X}_{\tau_{t}})^*
\mathcal{\bar{E}}_{s,t}^{-1}\sigma^{-1}\left( \tau_s,  \bar{X}_{\tau_{s}}\right)^*\nonumber\\&+\psi D^l_{u}%
\bar{X}_{\tau_{t}}\sigma^{\prime}(\tau_t, \bar{X}_{\tau_{t}})^* \mathcal{\bar
{E}}_{s,t}^{-1}\sigma^{-1}\left( \tau_s,  \bar{X}_{\tau_{s}}\right) ^*\nonumber\\
&-\psi\sigma(\tau_t, \bar{X}_{\tau_{t}})^*
\mathcal{\bar{E}}_{s,t}^{-1}D^l_u\mathcal{\bar{E}}_{s,t}
\mathcal{\bar{E}}_{s,t}^{-1}
\sigma^{-1}\left( \tau_s,  \bar{X}_{\tau_{s}}\right) ^*\mathbf{1}_{u\le \tau_s}\nonumber\\
&+\psi\sigma(\tau_u,\bar{X}_{\tau_{u}})^*\sigma^{-1}\left( \tau_s,  \bar{X}_{\tau_{s}}\right)
\mathcal{\bar{E}}_{s,t}^{-1}
D^l_u\sigma^{-1}\left( \tau_s,  \bar{X}_{\tau_{s}}\right) ^*.
\label{eq:ft}
\end{align}
Here $D^l_u\sigma^{-1}\left( \tau_s,  \bar{X}_{\tau_{s}}\right) ^*=\sum_{k=1}^dD^l_u\bar{X}_{\tau_{s}}^k\left(\sigma^{-1}\partial_{x_k}\sigma\sigma^{-1}\left( \tau_s,  \bar{X}_{\tau_{s}}\right)\right) ^*$.
One has then to get an upper bound for the $L^\rho$-norm of each term. As many of the arguments are repetitive, we show the reader only some of  the arguments that are involved.
Let us start with the first term in \eqref{eq:ft}. We have
\begin{align*}
D_u\psi&=\varphi^{\prime}\left(  |\mathcal{E}_{0,t}^{-1}\left(  \mathcal{E}%
_{0,t}-\mathcal{\bar{E}}_{0,t}\right) |^2 \right)  D_{u}\left[ | \mathcal{E}%
_{0,t}^{-1}\left(  \mathcal{E}_{0,t}-\mathcal{\bar{E}}_{0,t}\right)  |^2\right]
\end{align*}
and $D_{u}\left[ | \mathcal{E}%
_{0,t}^{-1}\left(  \mathcal{E}_{0,t}-\mathcal{\bar{E}}_{0,t}\right)  |^2\right]=-2\Tr\left[\left(\mathcal{E}%
_{0,t}^{-1}\left(  \mathcal{E}_{0,t}-\mathcal{\bar{E}}_{0,t}\right)\right)^* 
\left(\mathcal{E}%
_{0,t}^{-1}D_{u}\mathcal{E}%
_{0,t}
\mathcal{E}%
_{0,t}^{-1}\mathcal{\bar{E}}_{0,t}
-\mathcal{E}%
_{0,t}^{-1}D_u\mathcal{\bar{E}}_{0,t}\right)
\right]$ .
From the estimates in~\eqref{eq:propE} and~\eqref{eq:A13}, we obtain
\begin{equation}
\label{eq:Dest}
\sup_{u\in [0,t]}\left\Vert D_u\psi \right\Vert
_{\rho}\leq \|\varphi^{\prime}\|_\infty C(\rho).
\end{equation}
Note that  
 if $\varphi^{\prime}\left(  |\mathcal{E}_{0,t}^{-1}\left(  \mathcal{E}%
_{0,t}-\mathcal{\bar{E}}_{0,t}\right)  |^2\right) \not=0$ then $\psi\neq0$ and, reasoning like in \eqref{majosigxisig-1}, we have
\begin{align*}
&  \E \left[ \left |D_u\psi 
\sigma(\tau_t, \bar{X}_{\tau_{t}})^*
\mathcal{\bar{E}}_{s,t}^{-1} \sigma^{-1}\left( \tau_s,  \bar{X}_{\tau_{s}}\right)^*\right|^\rho\right] 
 \le  C
\left\Vert D_u\psi \right\Vert
_{2\rho}^\rho
  \E \left[ \left (\left| 
\mathcal{{E}}_{0,t}^{-1}\right |\left |\mathcal{\bar{E}}_{0,\eta (s)} \right|\right )^{2\rho}\right]^{1/2}.
\end{align*}
Similar bounds hold for
the three other terms. Note that the highest requirements on the derivatives
of~$b$ and~$\sigma$ will come from the terms involving $D_u\bar{\mathcal{E}}$
in (\ref{eq:ft}). Gathering all the upper bounds, we  get that using \eqref{upper_bound_1p} then $\left\Vert \psi\sigma(\tau_t, \bar{X}_{\tau_{t}})^*\mathcal{\bar{E}}_{\cdot,t}^{-1}
\sigma^{-1}\left (\bar
{X}_{\tau_{\cdot}}\right )^*\right\Vert _{1,\rho}^\rho
\le C(t^{\rho/2}+t^\rho) \le  Ct^{\rho/2}$ since $0\le t\le T$.
From~\eqref{controle_normep}, we finally obtain 
\begin{equation}
\label{eq:S}
\left\Vert \int_{0}^{t}\psi\sigma(\tau_u, \bar{X}_{\tau_{u}})^*
\mathcal{\bar{E}}_{s,t}^{-1}\sigma^{-1}\left(\tau_s,
\bar{X}_{\tau_{s}}\right)  ^*\delta W_{s}\right\Vert
_{\rho}\leq C(\rho)t^{1/2}.
\end{equation}

We are now in position to conclude. Using Jensen's inequality, the results~\eqref{espcond},~\eqref{eq:S},~\eqref{eq:Dest},~\eqref{vitesse_forte},~\eqref{eq:propE},  and the definition of
$\varphi$ together with Chebyshev's inequality, we have for any $k>0$ that there exists a constant $C\equiv C(k)$ such that
\begin{align*}
& \mathbb{E}\left[  \left\vert \mathbb{E}\left[  \left.  W_{t}-W_{\tau_{t}%
}\right|  \bar{X}_{t}\right]  \right\vert ^{\rho}\right]  
\\
& \leq C \left(  t^{-\rho}(t-\tau_{t})^{\rho}%
\left\Vert \int_{0}^{t}\psi\sigma(\tau_t, \bar{X}_{\tau_{t}}%
)^* \mathcal{\bar{E}}_{s,t}%
^{-1}\sigma^{-1}\left( \tau_s,  \bar{X}_{\tau_{s}}\right) ^*\delta W_{s}\right\Vert _{\rho}^\rho\right.  
+(t-\tau_{t})^{\rho-1}%
\int_{\tau_{t}}^{t}\left\Vert D_u\psi \right\Vert _{\rho}^\rho du 
\\
&\left. +\sqrt{\E(|W_t-W_{\tau_t}|^{2\rho})} \left(\E(|\mathcal{{E}}_{0,t}-\mathcal{\bar{E}}_{0,t}|^{2k})\E(|\mathcal{{E}}_{0,t}^{-1}|^{2k})\right)^{1/4} \right)\\
& \leq C\left(  t^{-\rho/2}(t-\tau
_{t})^{\rho}+(t-\tau
_{t})^{\rho}+\left(\frac{1}{N}\right)^{\left( 2\rho+k(1\wedge 2\gamma)\right)  /4}\right)  
 \leq C\left(  \frac{(t-\tau_t)^\rho}{t^{\rho/2}%
}+\frac{1}{N^\rho}+\frac{1}{N^{\frac{  \rho}{2}+\frac{k(1\wedge 2\gamma)}{4} }}\right)  .
\end{align*}
Taking $k$ big enough, the conclusion follows.
\end{adem}

\begin{adem}[of Lemma~\ref{lemme_majorations}]
 The finiteness of $\sup_{0\leq s\leq t \le T}{\mathbb{E}}\left[
 |\mathcal{E}_{s,t}|^{\rho}\right]+\sup_{0\leq s\leq t \le T}{\mathbb{E}}\left[
 |\mathcal{\bar{E}}_{s,t}|^{\rho}\right]$ is obvious since ${\nabla b}$ and $\sigma'$ are bounded. The upper bound for $\sup_{0\leq s\leq t \le T}{\mathbb{E}}\left[
| \mathcal{E}_{s,t}^{-1}|^\rho\right]$ is obtained using the same method of proof as in Theorem
48, Section V.9,  p.320  in \cite{Protter},  together with Gronwall's lemma.

 The estimate (\ref{eq:A13}) on $D_u{\mathcal{E}}$ is given, for example, by
 Theorem 2.2.1 in \cite{N} for time independent coefficients. The same method of proof works for our case. 
In fact, 
let us remark that $\mathcal{{E}}$ satisfies \eqref{eq:E} and that $\mathcal{\bar{E}}$ satisfies 
\begin{align*}
&\mathcal{\bar{E}}_{\eta_u,t}=I+\int_{\eta_u}^t\mathcal{\bar{E}}_{\eta_u,\tau_s}\sigma' (\tau_s,\bar{X}_{\tau_s})
dW_s+\int_{\eta_u}^t\mathcal{\bar{E}}_{\eta_u,\tau_s}
{\nabla b} (\tau_s, \bar{X}_{\tau_s})
ds.
\end{align*}
On the other hand, we have  for $\eta(s)\le u\le t$
\begin{align}
\nonumber
D^l_u\bar{\mathcal{E}}_{\eta_s,t}&=
\bar{\mathcal{E}}_{\eta_s,\tau_u}\sigma_l'(\tau_u, \bar{X}_{\tau_u})
+\int_{\eta_u}^t \left[ 
\bar{\mathcal{E}}_{\eta_s,\tau_r}D^l_u\sigma'(\tau_r, \bar{X}_{\tau_r})
+
D^l_u\bar{\mathcal{E}}_{\eta_s,\tau_r} \sigma'(\tau_r, \bar{X}_{\tau_r})\right] dW_r
\\
&+\int_{\eta_u}^t \left[ 
\bar{\mathcal{E}}_{\eta_s,\tau_r}D_u^l{\nabla b}(\tau_r, \bar{X}_{\tau_r})
+
D^l_u\bar{\mathcal{E}}_{\eta_s,\tau_r}{\nabla b}(\tau_r, \bar{X}_{\tau_r})\right]dr.
\label{eq:DbarE}
\end{align}
In order to obtain \eqref{eq:A13}, we
use~\eqref{eq:propE}, $b \in C^{\gamma,2}_b(\R^d)$,
$\sigma \in C^{\gamma,2}_b(\mathcal{M}_d(\R))$
and Gronwall's lemma. In fact, for example, one applies 
the $L^\rho(\Omega)$-norm to \eqref{eq:DbarE}, then using H\"older's inequality one obtains \eqref{eq:A13} if one uses the chain rule for stochastic derivatives, \eqref{eq:DbX} and \eqref{eq:propE}. Finally using 
$\mathcal{\bar{E}}_{u,t}=\mathcal{\bar{E}}_{\eta(u),t}$, one obtains  \eqref{eq:A13} for $\bar{\mathcal{E}}$.

Furthermore, \eqref{vitesse_forte} can be easily obtained by noticing that
$(\bar{X}_t,\bar{\mathcal{E}}_{0,t})$ is the Euler scheme for the SDE
$(X_t,\mathcal{E}_{0,t})$ which has  coefficients Lipschitz continuous in space and $\gamma$-H\"older continuous in time,  and
by using the strong convergence order of $\frac{1}{2}\wedge \gamma$ (see e.g.~Proposition~14 \cite{Faure}).

\end{adem}

\end{document}